\documentclass[smallcondensed,envcountsect]{svjour3}
\usepackage{latexsym, enumerate, amsmath, amsfonts, amscd, amssymb, verbatim}
\usepackage{graphicx,graphics}
\usepackage{xcolor}
\counterwithin{equation}{section}
\begin{document}

\newtheorem{Theorem}{Theorem}[section]
\newtheorem{Proposition}[Theorem]{Proposition}
\newtheorem{Remark}[Theorem]{Remark}
\newtheorem{Lemma}[Theorem]{Lemma}
\newtheorem{Corollary}[Theorem]{Corollary}
\newtheorem{Definition}[Theorem]{Definition}
\newtheorem{Example}[Theorem]{Example}
\renewcommand{\theequation}{\thesection.\arabic{equation}}
\normalsize

\setcounter{equation}{0}

\title{Properties of Generalized Polyhedral Convex Multifunctions}
\author{Nguyen Ngoc Luan \and Nguyen Mau Nam \and Nguyen Dong Yen}
\institute{N. N. Luan \at
Department of Mathematics and Informatics, Hanoi National University of Education, 136 Xuan Thuy, Hanoi, Vietnam. Institute of Mathematics, Vietnam Academy of Science and Technology, 18 Hoang Quoc Viet, Hanoi 10307.\\
luannn@hnue.edu.vn
\and
N. M. Nam \at
Fariborz Maseeh Department of Mathematics and Statistics, Portland State University, Portland, OR 97207, USA.\\ Research of this author was partly supported by the USA National Science Foundation under grant DMS-2136228.\\
mnn3@pdx.edu
\and
N. D. Yen \at
Institute of Mathematics, Vietnam Academy of Science and Technology, Hanoi, Vietnam\\
ndyen@math.ac.vn}
\date{Received: date / Accepted: date}

\maketitle
\medskip

\begin{abstract}
This paper presents a study of generalized polyhedral convexity under basic operations on multifunctions. We address the preservation of generalized polyhedral convexity under sums and compositions of multifunctions, the domains and ranges of generalized polyhedral convex multifunctions, and the direct and inverse images of sets under such mappings. Then we explore the class of optimal value functions defined by a generalized polyhedral convex objective function and a generalized polyhedral convex constrained mapping. The new results provide a framework for representing the relative interior of the graph of a generalized polyhedral convex multifunction in terms of the relative interiors of its domain and mapping values in locally convex topological vector spaces. Among the new results in this paper is a significant extension of a result by Bonnans and Shapiro on the domain of generalized polyhedral convex multifunctions from Banach spaces to locally convex topological vector spaces.
\end{abstract}


\section{Introduction}
The concept of polyhedral convex sets or convex polyhedra can be traced back to ancient Greece, where Plato discussed the five regular polyhedra in his book ``Timaeus." However, it wasn't until the 19th century that the study of convex polyhedra gained significant interest due to their crucial role in the theory of linear programming and their connection to convex analysis and optimization. The notion of polyhedral convex sets has since been used to define polyhedral convex functions and polyhedral convex multifunctions, which require their epigraphs and graphs, respectively, to be polyhedral convex sets. Polyhedral convex sets, functions, and multifunctions have many nice properties that can be used in convex analysis and optimization, making them valuable in many applications; see, e.g.,~\cite{Lee_Tam_Yen_2005,mordukhovich_nam_2021,Qui_Yen_NA_2011,Rockafellar_1970,Zheng_Yang_2008}.

\medskip
The important role of polyhedral convex sets in optimization and other areas has led to the development of a more general concept called {\em generalized polyhedral convex set}  in the framework of locally convex Hausdorff topological vector spaces. This concept was introduced by Bonnans and Shapiro in their book ``Perturbation Analysis of Optimization Problems'',  where they defined a generalized polyhedral convex set as the intersection of a polyhedral convex set and a closed affine subspace \cite[Definition~2.195]{Bonnans_Shapiro_2000}. These more general sets allow for a broader range of applications in optimization and other fields, particularly in the theories of generalized linear programming and quadratic programming~\cite[Sections~2.5.7 and 3.4.3]{Bonnans_Shapiro_2000}. 

It is well known that any infinite-dimensional normed space equipped with the weak topology is not metrizable, but it is a locally convex Hausdorff topological vector space. Similarly, the dual space of any infinite-dimensional normed space equipped with the weak$^*$ topology is not metrizable, but it is a locally convex Hausdorff topological vector space. The just mentioned two models provide us with the most typical examples of locally convex Hausdorff topological vector spaces, whose topologies cannot be given by norms. Therefore, it is necessary to study the generalized polyhedral convex set in locally convex Hausdorff topological vector spaces.

A generalized polyhedral convex function and a generalized polyhedral convex multifunction can be defined accordingly by requiring their epigraphs and graphs, respectively, to be generalized polyhedral convex sets. So, the concept of generalized polyhedral convex set has proved to be very useful in many issues of convex analysis and applications; see~\cite{Ban_Mordukhovich_Song_2011,Ban_Song_2016,Gfrerer_2013,Gfrerer_2014,Lim_Tuan_Yen_2022a,Lim_Tuan_Yen_2022b,Luan_Acta_2018,Luan_Nam_Thieu_Yen,Luan_Yen_Optim_2020,Yen_Yang_2016,Zheng_2009,Zheng_Ng_2014,Zheng_Yang_2021} and the references therein. The paper of Luan, Yao, and Yen \cite{Luan_Yao_Yen_2018} can be seen as a comprehensive study on generalized polyhedral convex sets, generalized polyhedral convex functions on locally convex Hausdorff topological vector spaces, and the related constructions such as sum of sets, sum of functions, directional derivative, infimal convolution, normal cone, conjugate function, subdifferential, sum rule. Some results of~\cite{Luan_Yao_Yen_2018} can be considered as adequate extensions of the corresponding classical results in \cite[Section~19]{Rockafellar_1970}.

\medskip
The present paper explores many properties of generalized polyhedral convex multifunctions with refinements to the case of polyhedral convex ones. We first examine the generalized polyhedral convexity of the domains and ranges as well as the direct and inverse images of generalized polyhedral convex sets under these mappings. Our new results include an  answer to the question of whether the domain of a generalized polyhedral convex multifunction is also a generalized polyhedral convex set in locally convex topological vector spaces. This is an important extension of a related result by Bonnans and Shapiro in the Banach space setting (see~\cite[Theorem~2.207]{Bonnans_Shapiro_2000}). Then we study the preservation of polyhedral convexity for multifunctions under various operations. We provide our answers to another question asking if the sum or composition of two generalized polyhedral convex multifunctions remains a generalized polyhedral convex multifunction. The question has not been fully answered in the literature, even in the case of polyhedral convex mappings. We also study the generalized polyhedral convexity of the optimal value function, which is important in parametric optimization. The specific features of generalized polyhedral convex sets allow us to obtain a representation for their generalized relative interiors, which we use to study the generalized relative interiors of graphs of generalized polyhedral convex multifunctions in locally convex topological vector spaces. Our developments have great potential applications to the theory of generalized differentiation involving generalized polyhedral convex sets, functions, multifunctions, and to optimization theory.

\medskip
The paper is structured as follows. Section 2 introduces basic notions and results related to generalized polyhedral convex sets and multifunctions. In Section 3, we discuss some properties of generalized polyhedral convex multifunctions including their domains and ranges, as well as the direct and inverse images of generalized polyhedral convex sets under such mappings. The stability of generalized polyhedral convexity under basic operations is presented in Section 4. Section 5 is devoted to the study of generalized relative interiors of generalized polyhedral convex sets. We obtain a representation for the relative interior of the graph of a generalized polyhedral convex multifunction in locally convex topological vector spaces.

\medskip
In the sequel,  $X$, $Y$, and $Z$ are assumed to be  locally convex Hausdorff topological vector spaces. We use the notation~$X^*$ to denote the topological dual space of $X$, and $\langle x^*, x \rangle$ to represent the value of $x^* \in X^*$ at $x \in X$. For a subset $\Omega \subset X$, we denote its topological closure and interior by $\overline{\Omega}$ and ${\rm int}\,\Omega$, respectively. The same notation is used for subsets of $X^*$. The cone (resp., the linear subspace) generated by a set $\Omega \subset X$ are denoted by ${\rm cone}\,\Omega$ (resp., ${\rm span}\,\Omega$).

\section{Preliminaries}\markboth{\centerline{\sc Preliminaries}}{\centerline{\sc N.~N.~Luan, N.~M.~Nam, and N.~D.~Yen}}
\setcounter{equation}{0}

This section recalls the definitions of generalized polyhedral convex sets, functions, and multifunctions, as well as some basic notations and results that will be used throughout the paper. The readers are referred to \cite{Bonnans_Shapiro_2000} for more details.

\medskip 
Let $X_0 \subset X$ be a closed linear subspace. Recall that the codimension of $X_0$ is the dimension of the quotient space~$X/X_0$ (see~\cite[p. 106]{Rudin_1991}). In the lemma below, we present two well-known characterizations of finite-codimensional linear subspaces and provide a detailed proof for the convenience of the readers. 
	
\begin{Lemma}\label{Lemma_codim}
Let $X_0 \subset X$ be a closed linear subspace. Then the following statements are equivalent:
	\begin{description}
	\item{\rm (a)} $X_0$ is finite-codimensional.
	\item{\rm (b)} There exists a finite-dimensional linear subspace~$X_1$ of~$X$ such that $X_0+X_1=X$ and $X_0 \cap  X_1= \{0\}$.
	\item{\rm (c)} There exists a continuous linear mapping $T$ from $X$ to a locally convex Hausdorff topological vector space $W$ such that $W$ is finite-dimensional and $X_0={\rm ker}\,T$.  
	\end{description}
\end{Lemma}
\noindent {\bf Proof.} ${\rm (a)} \Longrightarrow  {\rm (b)}$ See the proof of \cite[Lemma~4.21(b)]{Rudin_1991}.\\[1ex]
\noindent
${\rm (b)} \Longrightarrow  {\rm (c)}$ Let $\pi_0\colon X \rightarrow X/X_0$, $x \mapsto x+X_0$ for all $x \in X$, be the canonical projection from $X$ onto the quotient space $X/X_0$. Consider further the linear operator $\Phi_0\colon X/X_0 \rightarrow X_1$ defined as follows. For any $x\in X$, there is a unique representation $x=x_0+x_1$, where $x_0\in X_0$ and $x_1\in X_1$. Then we set $\Phi_0(x)=x_1$ and observe that $\Phi_0$ is bijective. On one hand, by \cite[Theorem~1.41(a)]{Rudin_1991}, $\pi_0$ is a continuous linear mapping. On the other hand, $\Phi_0$ is a homeomorphism by \cite[Lemma~2.5]{Luan_Yen_Optim_2020}. Thus, the operator $T\colon X \rightarrow X_1$ given by $T=\Phi_0\circ \pi_0$ is linear and continuous with  ${\rm ker}\,T=X_0$. The proof of this implication is complete by letting $W=X_1$ and observing that $W$ is a locally convex Hausdoff topological vector space that is finite-dimensional. \\[1ex]
\noindent
${\rm (c)} \Longrightarrow  {\rm (a)}$ It is clear that the operator $\Phi\colon X/X_0 \rightarrow T(X)$, $x +X_0 \mapsto T(x)$ for all $x \in X$, is a bijective linear mapping. Since $T(X)$ is a linear subspace of the finite-dimensional space $W$, one has ${\rm dim}\,T(X) < \infty$. Therefore, $X_0$ is finite-codimensional because ${\rm dim} (X/X_0)={\rm dim}\,T(X)<\infty$.
$\hfill\Box$

\medskip
A subset $D \subset X$ is said to be a {\em generalized polyhedral convex set}, or a \textit{generalized convex polyhedron}, if there exist $x^*_i \in X^*$, $\alpha_i \in \mathbb R$, $i=1,2,\ldots,m$, and a closed affine subspace $L \subset X$ such that
\begin{equation}\label{eq_def_gpcs}
	D=\big\{ x \in X \;\big|\; x \in L,\ \langle x^*_i, x \rangle \leq \alpha_i,\  i=1,\ldots,m\big\}.
\end{equation}
If $D$ can be represented in the form \eqref{eq_def_gpcs} with $L=X$, then we say that it is a \textit{polyhedral convex set}, or a \textit{convex polyhedron}.

\medskip
Let $D$ be given as in~\eqref{eq_def_gpcs}. By \cite[Remark~2.196]{Bonnans_Shapiro_2000}, there exists a continuous surjective linear mapping $A$ from $X$ to a locally convex Hausdorff topological vector space $Z$ and a vector $z \in Z$ such that $L=\big\{x \in X \mid A(x)=z  \big\}$. Then $D$ can be represented by
\begin{equation*}
	D=\big\{ x \in X \; \big|\; A(x)=z,\  \langle x^*_i, x \rangle \leq \alpha_i,\  i=1,\ldots,m\big\}.
\end{equation*}
If $D$ is a polyhedral convex set in $X$, then one can choose $Z=\{0\}$, $A\equiv 0$, and $z=0$.

\medskip
It follows from the definition that every generalized polyhedral convex set is a closed set. If~$X$ is finite-dimensional, a subset $D\subset X$ is a generalized polyhedral convex set if and only if it is a polyhedral convex set. In that case, we can represent a given affine subspace $L\subset X$ as the solution set of a system of finitely many linear inequalities.

\medskip
The following representation theorem for generalized convex polyhedral sets in the spirit of~\cite{Rockafellar_1970} is crucial for our subsequent proofs.

\begin{Theorem}\label{rep_gpcs}{\rm (See \cite[Theorem~2.7]{Luan_Yen_Optim_2020} and \cite[Lemma~2.12]{Luan_Yao_Yen_2018})} Suppose that $D$ is a nonempty subset of $X$. The set $D$ is generalized polyhedral convex (resp., polyhedral convex) if and only if there exist $u_1, \ldots, u_k \in X$, $v_1, \ldots, v_{\ell} \in X$, and a closed linear subspace (resp., a closed linear subspace of finite codimension) $X_0 \subset X$ such that
	\begin{eqnarray*}
		\begin{array}{ll}
			D=\Big\{ \sum\limits_{i=1}^k \lambda_i u_i + \sum\limits_{j=1}^\ell \mu_j v_j\; \big | \;  & \lambda_i \geq 0, \, \forall i=1,\ldots,k, \, \sum\limits_{i=1}^k \lambda_i=1,\\
			&\mu_j \geq 0,\, \forall j=1,\ldots,\ell   \Big\}+X_0.
		\end{array}
	\end{eqnarray*}
\end{Theorem}

Given a function $f\colon X \rightarrow \Bar{\mathbb{R}}=[-\infty, \infty]$, recall that the epigraph of $f$ is defined by
\begin{equation*}
	\mbox{\rm epi}\, f=\big\{(x, \lambda)\in X\times \mathbb{R}\; \big |\; f(x)\leq \lambda\big\}.
\end{equation*}
The function $f$ is said to be a {\em lower semicontinuous} if $\mbox{\rm epi}\, f$ is a closed set in $X\times \mathbb{R}$. We also say that $f$ is a {\em generalized polyhedral convex function} (resp., a {\em polyhedral convex function}) if $\mbox{\rm epi}\, f$ is a generalized polyhedral convex set (resp., a polyhedral convex set) in $X\times \mathbb{R}$.

\medskip
Let $F\colon  X \rightrightarrows Y$ be a multifunction. The \textit{domain},  \textit{range}, and \textit{graph} of $F$ are defined, respectively, by $${\rm dom}\,F=\{x \in X \mid F(x) \neq \emptyset \}, \ \; {\rm rge}\,F=\{y \in Y \mid \exists x \in X  \ \, \textrm{such that} \ \, y \in F(x) \}$$ and $${\rm gph}\,F=\big\{(x, y) \in X \times Y \; \big| \; x \in {\rm dom}\,F,\, y \in F(x)\big\}.$$
It is clear that ${\rm dom}\,F=\pi_X({\rm gph}\,F)$, where $\pi_X\colon \, X \times Y \rightarrow X$ with $\pi_X(x,y)=x$ is the projection mapping from $X \times Y$ onto~$X$. Similarly, ${\rm rge}\,F=\pi_Y({\rm gph}\,F)$, where $\pi_Y\colon \, X \times Y \rightarrow Y$  with $\pi_Y(x,y)=y$ is the projection mapping from $X \times Y$ onto~$Y$. Observe that ${\rm  rge}\,F^{-1}={\rm dom}\,F$ and ${\rm  rge}\,F={\rm dom}\,F^{-1}$, where the inverse multifunction $F^{-1}\colon Y \rightrightarrows X$ of $F$ is given by $F^{-1}(y)=\{x \in X \mid y \in F(x)\}$.

\medskip
A multifunction $F\colon X \rightrightarrows Y$ is said to be  {\em generalized polyhedral} (resp., {\em polyhedral}) if ${\rm gph}\,F$ is a union of finitely many generalized polyhedral convex sets (resp., polyhedral convex sets) in $X \times Y$. If ${\rm gph}\,F$ is generalized polyhedral convex (resp., polyhedral convex), then we say that $F$ is \textit{generalized polyhedral convex} (resp., \textit{polyhedral convex}).

\medskip
Clearly, if $F\colon X \rightrightarrows Y$ is generalized polyhedral convex (resp., polyhedral convex), then $F^{-1}$ is also generalized polyhedral convex (resp., polyhedral convex).

\medskip
Let $F\colon X \rightrightarrows Y$ be a generalized polyhedral convex multifunction. If $A\colon X\times Y \rightarrow Z$ is a continuous linear mapping, then the formula $A_1(x)=A(x,0)$ for $x\in X$  (resp., the formula $A_2(y)=A(0,y)$ for $y\in Y$) defines a continuous linear mapping from $X$ to $Z$ (resp., a continuous linear mapping from $Y$ to $Z$), and one has
\begin{equation*}
	A(x,y)=A_1(x)+A_2(y), \ \;  (x,y)\in X\times Y.
\end{equation*}
Note also that $(X\times Y)^*=X^*\times Y^*$. Therefore, the graph of $F$ can be given by the formula
\begin{eqnarray}\label{def_Graph_F_1}
	\begin{array}{ll}
		{\rm gph}\,F=\big\{(x, y) \in X \times Y \; \big |\; &A_1(x)+A_2(y)=z,\\
		&\langle x^*_i, x \rangle + \langle y^*_i, y \rangle \leq \beta_i,\, i=1,\ldots,m\big\},
	\end{array}
\end{eqnarray}	
where $A_1$ (resp., $A_2$) is a continuous linear mapping from $X$ (resp., from $Y$) to $Z$, $z \in Z$, $x^*_i \in X^*$, $y^*_i \in Y^*$, $\beta_i \in \mathbb{R}$, for $i=1,\ldots,m$. Conversely, if the graph of a multifunction $F\colon X\rightrightarrows  Y$ can be represented by \eqref{def_Graph_F_1}, then $F$ is a generalized polyhedral convex multifunction. Note that if $\mbox{\rm gph}\, F$ is the emptyset, then by definition $F$ is a generalized polyhedral convex multifunction.  Clearly, if $F$ is a polyhedral convex multifunction, then one can choose $Z=\{0\}$, $A_1\equiv 0$, $A_2 \equiv 0$ and $z=0$. If the graph of a multifunction $F\colon X\rightrightarrows  Y$ can be represented by \eqref{def_Graph_F_1}, where $Z=\{0\}$, $A_1\equiv 0$, $A_2 \equiv 0$ and $z=0$, then $F$ is a polyhedral convex multifunction. 

\medskip 
We say that a continuous linear mapping $A\colon Y \rightarrow Z$ is \textit{closed under finite-codimensional subspaces} if $A(Y_0)$ is closed for every finite-codimensional closed linear subspace $Y_0 \subset Y$. Clearly, if $A$ is closed under finite-codimensional subspaces, then $A(Y)$ is closed. The converse is true if both $Y$ and $Z$ are Fr\'echet spaces (see Theorem~\ref{codimen_property} below). In particular, any continuous linear mapping between Banach spaces with a closed range is closed under finite-codimensional subspaces. This fact has been established by one of the two referees with direct proof. 

\section{Properties of Generalized Polyhedral Convex Multifunctions}\markboth{\centerline{\sc Properties of Generalized Polyhedral Convex Multifunctions}}{\centerline{\sc N.~N.~Luan, N.~M.~Nam, and N.~D.~Yen}}
\setcounter{equation}{0}

In this section, we study the domains and ranges of generalized polyhedral convex multifunctions, as well as the direct and inverse images of generalized polyhedral convex sets under such mappings. We will also discuss refinements of the obtained results in the case of polyhedral convex sets and multifunctions.

\medskip
First, we will extend a part of Theorem~2.207 in the book by Bonnans and Shapiro~\cite{Bonnans_Shapiro_2000}, which was given in a Banach space setting, to the case of generalized polyhedral convex multifunctions in locally convex Hausdorff topological vector spaces.

\begin{Theorem}\label{Property_domF} If the graph of a multifunction $F$ is described by~\eqref{def_Graph_F_1} in which the mapping $A_2$ is closed under finite-codimensional subspaces, then ${\rm dom}\,F$ is a generalized polyhedral convex set.
\end{Theorem}
\noindent {\bf Proof.} Without loss of generality we can assume that ${\rm gph}\,F$ is nonempty. Fix an element $(\bar x, \bar y) \in {\rm gph}\,F$. We observe that ${\rm gph}\,F=(\bar x, \bar y) + Q$ with
\begin{eqnarray*}
	\begin{array}{ll}
		Q=\big\{(x, y) \in X \times Y \; \big| \; & A_1(x)+A_2(y)=0, \\
		&\langle x^*_i, x \rangle + \langle y^*_i, y \rangle \leq \bar{\beta}_i,\, i=1,\ldots,m\big\},
	\end{array}
\end{eqnarray*}	
where $\bar{\beta}_i =\beta_i- \langle x^*_i, \bar{x} \rangle - \langle y^*_i, \bar{y} \rangle $ for $i=1,\ldots,m$.
Let $$W=\big\{(x, y) \in X \times Y \mid A_1(x)+A_2(y)=0\big\}.$$
Since $W_0=\big\{(x, y) \in W \mid \langle x^*_i, x \rangle + \langle y^*_i, y \rangle =0,\, i=1,\ldots,m\big\}$
is a closed linear subspace of finite codimension in $W$, one can find a finite-dimensional linear subspace~$W_1$ of~$W$ such that $W_0+W_1=W$ and $W_0 \cap  W_1= \{0\}$ by Lemma~\ref{Lemma_codim}. The subspace $W_1$ is closed by~\cite[Theorem 1.21(b)]{Rudin_1991}. Obviously,
\begin{equation*}
	Q_1=\big\{(x, y) \in W_1 \; \big|\; \langle x^*_i, x \rangle + \langle y^*_i, y \rangle \leq \bar{\beta}_i,\, i=1,\ldots,m\big\}
\end{equation*}	
is a polyhedral convex set in $W_1$. It is clear that $W_0+Q_1\subset Q$. The reverse inclusion is also true. Indeed, for each $(x, y) \in Q$ there exist $(x_0,y_0) \in W_0$ and $(x_1,y_1) \in W_1$ satisfying $(x, y)=(x_0,y_0)+(x_1,y_1)$. Since
\begin{equation*}
	\begin{aligned}
		\langle x^*_i, x_1 \rangle + \langle y^*_i, y_1 \rangle =&[\langle x^*_i, x \rangle + \langle y^*_i, y \rangle]-[\langle x^*_i, x_0 \rangle + \langle y^*_i, y_0 \rangle] \leq \bar{\beta}_i
	\end{aligned}
\end{equation*}
for all $i=1,\ldots,m$, one has $(x_1,y_1) \in Q_1$, so $(x, y) \in W_0+Q_1$. We have thus proved that $Q=W_0+Q_1$. Since $Q_1$ is a polyhedral convex set of the finite-dimensional space~$W_{1}$, invoking~\cite[Theorem 19.1]{Rockafellar_1970}, one can find $(x_1, y_1), \ldots, (x_k, y_k)$ in $Q_1$, $(u_1, v_1), \ldots, (u_{\ell}, v_{\ell})$ in~$W_1$ such that
\begin{eqnarray*}\label{rep_Q1}
	\begin{array}{ll}
		Q_1=\Big\{ \sum\limits_{i=1}^k \lambda_i (x_i, y_i) + \sum\limits_{j=1}^\ell \mu_j (u_j, v_j) \; \big|\; & \lambda_i \geq 0, \, \forall i=1,\ldots,k, \, \sum\limits_{i=1}^k \lambda_i=1, \\
		&\mu_j \geq 0,\, \forall j=1,\ldots,\ell \Big\}.
	\end{array}
\end{eqnarray*}
From what has already been said we obtain
\begin{eqnarray*}\label{rep_Graph_F_1}
	\begin{array}{ll}
		{\rm gph}\,F\!=\!(\bar x, \bar y)+\Big\{ \sum\limits_{i=1}^k \lambda_i (x_i, y_i) + \sum\limits_{j=1}^\ell \mu_j (u_j, v_j) \; \big|\; & \lambda_i \geq 0, \, \forall i=1,\ldots,k,  \,\sum\limits_{i=1}^k \lambda_i=1,  \\
		&\mu_j \geq 0,\, \forall j=1,\ldots,\ell   \Big\}\!+\!W_0.
	\end{array}
\end{eqnarray*}
Consequently,
\begin{eqnarray}\label{rep_dom_F_1}
	\begin{array}{ll}
		{\rm dom}\,F=\Big\{ \sum\limits_{i=1}^k \lambda_i (\bar{x}+x_i) + \sum\limits_{j=1}^\ell \mu_j u_j \; \big|\; &\lambda_i \geq 0, \, \forall i=1,\ldots,k, \sum\limits_{i=1}^k \lambda_i=1, \\
		&\mu_j \geq 0,\, \forall j=1,\ldots,\ell   \Big\}+\pi_X(W_0).
	\end{array}
\end{eqnarray}
Let $\widetilde{A}_1\colon X \to Z \times \mathbb{R}^m$, $\widetilde{A}_2\colon Y \to Z \times \mathbb{R}^m$  be continuous linear mappings defined, respectively, by
$$\widetilde{A}_1(x)=\big(A_1(x),\, \langle x^*_1, x \rangle, \ldots,  \langle x^*_m, x \rangle\big), \ \;  \widetilde{A}_2(y)=\big(A_2(y),\, \langle y^*_1, y \rangle, \ldots,  \langle y^*_m, y \rangle\big).$$ It follows that
\begin{eqnarray}\label{rep_PXW0}
	\begin{array}{ll}
		\pi_X(W_0)&=\big\{x \in X \;\big|\; \text{ there exists }\, y \in Y \, \, \text{such that}\, \, \widetilde{A}_1(x)+\widetilde{A}_2(y)=0\big\}\\
		&=\big\{x \in X \; \big|\; \widetilde{A}_1(x) \in \widetilde{A}_2(Y)\big\}\\
		&=\widetilde{A}_1^{-1}(\widetilde{A}_2(Y)).\\
	\end{array}
\end{eqnarray}

\noindent
Next, we will claim that the linear subspace $\widetilde{A}_2(Y)$ is closed in $Z \times \mathbb{R}^m$. Indeed, since $$Y_0=\big\{y \in Y \mid \langle y^*_i, y \rangle =0\ \, \mbox{\rm for all}\ \, i=1,\ldots,m\big\}$$
is a closed linear subspace of finite codimension in $Y$, one can find a finite-dimensional linear subspace~$Y_1$ of~$Y$ such that $Y_0+Y_1=Y$ and $Y_0 \cap  Y_1= \{0\}$ by Lemma~\ref{Lemma_codim}. The linear subspace $Y_1$ is closed by~\cite[Theorem 1.21(b)]{Rudin_1991}. We have
\begin{equation*}
\begin{array}{rl}
	\widetilde{A}_2(Y)=\widetilde{A}_2(Y_0+Y_1) & =\widetilde{A}_2(Y_0)+\widetilde{A}_2(Y_1)\\
	&=\big(A_2(Y_0) \times \{0\}\big)+\widetilde{A}_2(Y_1).
\end{array}
\end{equation*}
Since the mapping $A_2$ is closed under finite-codimensional subspaces, we see that $A_2(Y_0)$ is closed in $Z$ and hence $A_2(Y_0) \times \{0\}$ is a closed linear subspace of $Z \times \mathbb{R}^m$. As $\widetilde{A}_2(Y_1)$ is a finite-dimensional subspace of $Z \times \mathbb{R}^m$, by~\cite[Theorem~1.42]{Rudin_1991} we can assert that $\big(A_2(Y_0) \times \{0\}\big)+\widetilde{A}_2(Y_1)$ is closed. Thus, $\widetilde{A}_2(Y)$ is a closed linear subspace of $Z \times \mathbb{R}^m$.

Combining the result in the claim above with the continuity of the linear mapping~$\widetilde{A}_1$, we can assert that $\widetilde{A}_1^{-1}(\widetilde{A}_2(Y))$ is closed. Therefore, the subspace $\pi_X(W_0)$ is closed by formula~\eqref{rep_PXW0}. From~\eqref{rep_dom_F_1}, we conclude that ${\rm dom}\,F$ is a generalized polyhedral convex set by Theorem~\ref{rep_gpcs}.
$\hfill\Box$

\medskip
To derive the following result, we use a similar argument as in the proof of Theorem~\ref{Property_domF}.

\begin{Theorem}\label{Property_rgeF} If the graph of a multifunction $F$ is described by~\eqref{def_Graph_F_1} in which the mapping $A_1$ is closed under finite-codimensional subspaces, then ${\rm rge}\,F$ is a generalized polyhedral convex set.
\end{Theorem}

We can explore an interesting question related to Theorem~\ref{Property_domF}: \textit{Can the assumption that the mapping $A_2$ is closed under finite-codimensional subspaces be removed from this theorem?} To answer this question, let us provide an example.

\begin{Example}\label{Example1}{\rm Let $X=Y=C[a, b]$, $a<b$, be the linear space of continuous real-valued functions on the interval $[a, b]$ with the norm defined by $\|x\|=\max\limits_{t \in [a,b]}|x(t)|$. Let the continuous linear mappings $A_1\colon X \rightarrow X$ and $A_2\colon Y \rightarrow X$  be defined, respectively, by $A_1(x)=x$ and $\big(A_2(y)\big)(t)=\displaystyle \int_{a}^{t}y(\tau) d \tau, $ where integral is understood in the Riemannian sense. It is clear that 	
\begin{eqnarray}\label{Omega_ex1}\Omega=\big\{(x, y) \in X \times Y \mid A_1(x) + A_2(y)=0\big\}\end{eqnarray} is a closed linear subspace of $X \times Y$; hence it is a generalized polyhedral convex set in $X \times Y$. Let $F\colon X \rightrightarrows Y$ be the generalized polyhedral convex multifunction with ${\rm gph}\,F=\Omega$. Then we have
		\begin{eqnarray*}
			\begin{array}{ll}
				{\rm dom}\,F&=\pi_X({\rm gph}\,F)\\
				&=\big\{x \in X \; \big | \; \text{ there exists }\, y \in Y \, \, \text{such that}\, \, A_1(x) + A_2(y)=0\big\}\\
				&=\big\{x \in X \; \big | \; A_1(x) \in A_2(Y)\big\}\\
				&=\big\{x \in X \; \big | \; x \in A_2(Y)\big\}=A_2(Y).\\
			\end{array}
		\end{eqnarray*}
Since $A_2(Y)=\big\{ x \in C[a,b] \; \big | \; x \text{ is continuously differentiable on } (a,b),\, x(a)=0 \big\}$ is a non-closed linear subspace of $X$ (see \cite[Example 2.1]{Luan_APAN_2019} for details), ${\rm dom}\,F$ is not a generalized polyhedral convex set.}
\end{Example}

The following theorem addresses a special case of Theorem~\ref{Property_domF} in which $F$ is a polyhedral convex multifunction.

\begin{Theorem}\label{pcm}
If $F$ is a polyhedral convex multifunction, then ${\rm dom}\, F$ and ${\rm rge}\, F$ are polyhedral convex sets.
\end{Theorem}
{\bf Proof.} We can assume that ${\rm gph}\,F$ is given by \eqref{def_Graph_F_1} with $Z=\{0\}$, $A_1\equiv 0$, $A_2 \equiv 0$ and $z=0$. Arguing similarly as in the proof of Theorem~\ref{Property_domF}, to prove that ${\rm dom}\, F$ is a polyhedral convex set in $X$, we need only to show that $\pi_X(W_0)$ is a closed linear subspace of finite codimension of $X$. In the notation of the proof of Theorem~\ref{Property_domF}, observe that $\widetilde{A}_2 (Y)$ is a linear subspace of the finite dimensional space $\{0\} \times \mathbb{R}^m$. Therefore, $\widetilde{A}_1^{-1}(\widetilde{A}_2(Y))$ is a closed linear subspace of finite codimension in $X$. Thus, from~\eqref{rep_PXW0} it follows that $\pi_X(W_0)$ is a closed linear subspace of finite codimension in $X$.

The fact that ${\rm rge}\, F$ a polyhedral convex set in $Y$ is obtained from the above result by considering the multifunction $F^{-1}$, which is polyhedral convex by the assumption of the theorem, and applying the formula ${\rm  rge}\,F={\rm dom}\,F^{-1}$. $\hfill\Box$

\medskip
As a direct consequence of Theorem~\ref{pcm}, we now present a property of the projection of a polyhedral convex set in a product space onto each component of the latter.

\begin{Corollary}\label{pcs_cr}
If $P$ is a polyhedral convex set in $X \times Y$, then $\pi_X(P)$ is a polyhedral convex set in $X$ and $\pi_Y(P)$ is a polyhedral convex set in $Y$.
\end{Corollary}	
{\bf Proof.} Suppose that $P\subset X \times Y$ is a polyhedral convex set. Let $G\colon X \rightrightarrows Y$ be the multifunction defined by
$$G(x)=\big\{y \in Y \, \big | \, (x,y) \in P\big\}, \ \; x \in X.$$
Since ${\rm gph}\,G=P$, we see that $G$ is a polyhedral convex multifunction. As $\pi_X(P)={\rm dom}\,G$ and $\pi_Y(P)={\rm rge}\,G$, the desired properties follow from Theorem~\ref{pcm}.
$\hfill\Box$

\medskip
The proposition below gives us a useful result concerning the generalized polyhedral convexity of the function values of a generalized polyhedral convex multifunction.

\begin{Proposition}\label{gpcs_image1} If $F\colon X \rightrightarrows Y$ is a generalized polyhedral convex multifunction, then $F(x)$ is a generalized polyhedral convex set in $Y$ for every $x \in X$.
\end{Proposition}
{\bf Proof.} We can assume that the graph of $F$ can be given by~\eqref{def_Graph_F_1}. Taking any element $\bar x \in X$, we have
\begin{eqnarray*}
	\begin{array}{ll}
		F(\bar x)&=\big\{y \in Y \; \big |\; (\bar x,y) \in {\rm gph}\,F\big\}\\
		 &=\big\{y \in Y  \; \big |\;A_1(\bar x)+A_2(y)=z,\, \langle x^*_i, \bar x \rangle + \langle y^*_i, y \rangle \leq \beta_i,\, i=1,\ldots,m\big\}\\
		 &=\big\{y \in Y  \; \big |\;A_2(y)=z-A_1(\bar x),\,   \langle y^*_i, y \rangle \leq \beta_i-\langle x^*_i, \bar x \rangle,\, i=1,\ldots,m\big\}.		 
	\end{array}
\end{eqnarray*}
This clearly implies that $F(\bar x)$ is a generalized polyhedral convex set in $Y$. 
$\hfill\Box$

\medskip
We now show that the image of a generalized polyhedral convex set under a polyhedral convex multifunction is  a polyhedral convex set.

\begin{Proposition}\label{pcs_image1} Let $F\colon X \rightrightarrows Y$ be a polyhedral convex multifunction. If $C \subset X$ is a generalized polyhedral convex set in $X$, then $F(C)$ is a polyhedral convex set in $Y$. In particular, if $P\subset X$ is a polyhedral convex set in~$X$, then $F(P)$ is a polyhedral convex set in $Y$.
\end{Proposition}
{\bf Proof.} Without loss of generality, suppose that both sets ${\rm gph}\,F$ and $C$ are nonempty.
By the assumptions, we can assume that $${\rm gph}\,F=\big\{(x, y) \in X \times Y \, \; \big | \;
\, \langle x^*_i, x \rangle + \langle y^*_i, y \rangle \leq \alpha_i,\, i=1,\ldots,p\big\}$$
and
$C=\big\{ x \in X \; \big | \; A(x)=z,\  \langle u^*_j, x \rangle \leq \beta_j,\  j=1,\ldots,q\big\},$ where $x^*_i \in X^*, y^*_i \in Y^*, \alpha_i \in \mathbb{R}$ for $i=1,\ldots,p$, $A\colon X \to Z$ is a continuous linear mapping, $z \in Z$, $u^*_j \in X^*, \beta_j \in \mathbb{R}$ for $j=1,\ldots,q$. Set $X_0={\rm ker}\,A$, $\Omega={\rm gph}\,F \cap (C \times Y)$, and select a point $\bar{x} \in C$. Let
\begin{eqnarray}\label{Omega0}
	\begin{array}{ll}
		\Omega_0=\big\{(x_0, y) \in X_0 \times Y \; \big | \; &
		\langle x^*_{0,i}, x_0 \rangle + \langle y^*_{i}, y \rangle \leq \alpha_{0,i},\, i=1,\ldots,p\\
		&\langle u^*_{0,j}, x_0 \rangle \leq \beta_{0,j},\  j=1,\ldots,q \big\},
	\end{array}
\end{eqnarray}
where $x^*_{0,i}$ denotes the restriction of $x^*_i$ to $X_0$, $\alpha_{0,i}=\alpha_i-\langle x^*_i, \bar{x} \rangle$ for $i=1,\ldots,p$,  $u^*_{0,j}$ denotes the restriction  of $u^*_j$ to $X_0$, and $\beta_{0,j}=\beta_j-\langle u^*_j, \bar{x} \rangle$ for  $j=1,\ldots,q$. It is easily verified that
$\Omega=(\bar x,0)+\Omega_0.$ Since $F(C)=\pi_Y(\Omega)$, this implies that
\begin{eqnarray}\label{F(C)} F(C)=\pi_Y(\bar x,0)+\pi_Y(\Omega_0)=\pi_Y(\Omega_0).
\end{eqnarray}
By~\eqref{Omega0}, the set $\Omega_0$ is  polyhedral convex in $X_0 \times Y$. According to Corollary~\ref{pcs_cr}, $\pi_Y(\Omega_0)$ is a polyhedral convex set in $Y$. Hence, from~\eqref{F(C)} it follows that $F(C)$ is a polyhedral convex set in $Y$.
$\hfill\Box$

The next example shows that both assertions of Proposition~\ref{pcs_image1} are false if $F$ is merely a generalized polyhedral convex multifunction.

\begin{Example}\label{gtg} {\rm Let $X, Y,$ and $F$ be the same as in Example~\ref{Example1}. Since $F$ is a generalized polyhedral convex multifunction, its inverse  $F^{-1}\colon Y \rightrightarrows X$ is also a generalized polyhedral convex multifunction. Obviously, $Y$ is a polyhedral convex set in $Y$ and $$F^{-1}(Y)={\rm rge}\,F^{-1}={\rm dom}\,F.$$ Since ${\rm dom}\,F$ is not a generalized polyhedral convex set, the image of $Y$ via the generalized polyhedral convex multifunction $F^{-1}$ is not a generalized polyhedral convex set.}
\end{Example}

As a direct consequence of Proposition~\ref{pcs_image1}, the corollary below addresses the inverse image of a generalized polyhedral convex set under a polyhedral convex multifunction.

\begin{Corollary} Let $F\colon X \rightrightarrows Y$ be a polyhedral convex multifunction. If $D \subset Y$ is a generalized polyhedral convex set, then $F^{-1}(D)$ is a polyhedral convex set in $X$. In particular, if $Q\subset Y$ is a polyhedral convex set, then $F^{-1}(Q)$ is a polyhedral convex set in  $X$.
\end{Corollary}
{\bf Proof.} By the assumptions made, $F^{-1}\colon Y \rightrightarrows X$ is a polyhedral convex multifunction and $D$ is a generalized polyhedral convex set in $Y$. Then, by Proposition~\ref{pcs_image1} we can assert that $F^{-1}(D)$ is a polyhedral convex set in $X$.
$\hfill\Box$

\medskip
Example~\ref{gtg} has demonstrated that the conclusions of Proposition~\ref{pcs_image1} do not hold if $F$ is just assumed to be a generalized polyhedral convex multifunction. Nevertheless, when $F$ is \textit{a surjective continuous linear mapping between Fr\'echet spaces} (a specific type of generalized polyhedral convex multifunction), we have the following intriguing result. Recall~\cite[p. 9]{Rudin_1991} that a locally convex Hausdorff topological vector space is said to be a \textit{Fr\'echet space} if its topology $\tau$ is induced by a complete invariant metric $d$.

\begin{Theorem}\label{image_pcs} Let $X$ and $Y$ be Fr\'echet spaces and $T\colon X \rightarrow Y$ a surjective continuous linear mapping. If $D\subset X$ a polyhedral convex set, then $T(D)$ is a polyhedral convex set in~$Y$.
\end{Theorem}
{\bf Proof.} We can assume that the polyhedral convex set $D$ is nonempty. Then, according to Theorem~\ref{rep_gpcs}, there exist $u_1, \ldots, u_k \in X$, $v_1, \ldots, v_{\ell} \in X$, and a closed linear subspace of finite codimension $X_0$ in $X$ such that
\begin{eqnarray*}
	\begin{array}{ll}
		D=\Big\{ \sum\limits_{i=1}^k \lambda_i u_i + \sum\limits_{j=1}^\ell \mu_j v_j \; \big| \; & \lambda_i \geq 0, \, \forall i=1,\ldots,k, \, \sum\limits_{i=1}^k \lambda_i=1,   \\
		&\mu_j \geq 0,\, \forall j=1,\ldots,\ell   \Big\}+X_0.
	\end{array}
\end{eqnarray*}
Thus, by the linearity of $T$ one has
\begin{eqnarray}\label{T_eq_rep_gpcs}
	\begin{array}{ll}
		T(D)=\Big\{ \sum\limits_{i=1}^k \lambda_i T(u_i) + \sum\limits_{j=1}^\ell \mu_j T(v_j) \; \big|\; & \lambda_i \geq 0, \, \forall i=1,\ldots,k,  \, \sum\limits_{i=1}^k \lambda_i=1,  \\
		&\mu_j \geq 0,\, \forall j=1,\ldots,\ell   \Big\}+T(X_0).
	\end{array}
\end{eqnarray}
Consider the quotient space $X/{\rm ker}\,T$ and the quotient mapping $\pi\colon  X\to X/{\rm ker}\,T$ (see, e.g.,~\cite[pp.~30-31]{Rudin_1991}) given by
\begin{equation*}
	\pi(x)=[x]=x+{\rm ker}\,T, \ \; x\in X.
\end{equation*}
According to~\cite[Theorem 1.41(a)]{Rudin_1991}, $\pi$ is linear, continuous, and open.

Let the linear mapping $\widetilde{T}\colon X/{\rm ker}\,T \to Y$ be defined by $\widetilde{T}([x])=T(x)$ for  $x\in X$. Since $T$ is a surjective continuous linear  mapping, $\widetilde{T}$ is a bijective continuous linear mapping. As $X$ is a Fr\'{e}chet space, by~\cite[Theorem 1.41(d)]{Rudin_1991} we know that $X/{\rm ker}\,T$ is also a Fr\'{e}chet space. So, thanks to Corollary~2.12(b) in~\cite{Rudin_1991},  ${\widetilde{T}}^{-1}\colon  Y \to X/{\rm ker}\,T$ is a continuous linear mapping. Since ${\rm ker}\,T \subset X$ is a closed linear subspace and $X_0 \subset X$ is a closed linear subspace of finite codimension, by~\cite[Lemma~2.15]{Luan_Yao_Yen_2018} we can infer that $X_0+{\rm ker}\,T$ is a closed linear subspace of $X$. Then, by the openness of $\pi$, the set $\pi\Big(X\setminus \big(X_0+{\rm ker}\,T\big)\Big)$ is open in $X/{\rm ker}\,T$. Consequently, the obvious equality
$$\pi\big(X_0+{\rm ker}\,T\big)=(X/{\rm ker}\,T)\setminus \pi\Big(X\setminus \big(X_0+{\rm ker}\,T\big)\Big)$$ shows that $\pi\big(X_0+{\rm ker}\,T\big)$ is a closed linear subspace of the quotient space. Since $$T(X_0)=({\widetilde{T}}^{-1})^{-1}\Big(\pi\big(X_0+{\rm ker}\,T\big)\Big),$$ the latter fact implies that $T(X_0)$ is a closed linear subspace of $Y$. As $\mbox{\rm codim}\,X_0 < \infty$ and $T$ is a surjective mapping, ${\rm codim}(T(X_0)) < \infty$. Hence $T(X_0)\subset Y$ is a closed linear subspace of finite codimension. Therefore, by the representation~\eqref{T_eq_rep_gpcs} and Theorem~\ref{rep_gpcs} we can assert that $T(D)$ is a polyhedral convex set in~$Y$.
$\hfill\Box$

\medskip
The last theorem of this section establishes the relationship between the closedness under finite-codimen\-sional subspaces and the closed range property of a continuous linear mapping between Fr\'echet spaces.

\begin{Theorem}\label{codimen_property} Let $X$ and $Y$ be Fr\'echet spaces and $T\colon X \to Y$ a continuous linear mapping. Then $T$ is closed under finite-codimensional  subspaces if and only if $T(X)$ is closed.  
\end{Theorem}
{\bf Proof.} It suffices to prove that if $T(X)$ is closed then $T$ is closed under finite-codimensional subspaces because the converse implication is obvious. Suppose that $X_0$ is a finite-codimensional closed linear subspace of $X$.  Using Lemma~\ref{Lemma_codim}, we have the representation $X=X_0+X_1$, where $X_0 \cap  X_1= \{0\}$ and~$X_1$ is a finite-dimensional subspace of $X$. By Theorem~\ref{rep_gpcs}, the linear subspace $X_0$ is a polyhedral convex set in $X$. Since $T(X)$ is closed, the continuous linear mapping $T \colon X \rightarrow T(X)$ is surjective between two Fr\'echet spaces. Thus, it follows from Theorem~\ref{image_pcs} that $T(X_0)$ is a polyhedral convex set in $T(X)$, so $T(X_0)$ is closed in $T(X)$. Since $T(X)$ is closed in $Y$, we see that $T(X_0)$ is closed in $Y$.  Therefore, $T$ is closed under finite-codimensional subspaces. $\hfill\Box$

\medskip 
Combining Theorem~\ref{codimen_property} with Theorem~\ref{Property_domF}, we get the following result. 
\begin{Corollary} If $Y$ and $Z$ are Fr\'echet spaces and the graph of a multifunction $F$ is given by~\eqref{def_Graph_F_1} in which the mapping $A_2$ has a closed range, then ${\rm dom}\,F$ is a generalized polyhedral convex set.
\end{Corollary}

\section{Generalized Polyhedral Convexity under Basic Operations}\markboth{\centerline{\sc Generalized Polyhedral Convexity under Basic Operations}}{\centerline{\sc N.~N.~Luan, N.~M.~Nam, and N.~D.~Yen}}
\setcounter{equation}{0}

In this section, we will examine the generalized polyhedral convex property in relation to basic operations on multifunctions. Our goal is to study the preservation of the generalized polyhedral convexity under sums and compositions of multifunctions. We will also study an important class of extended-real-valued functions known as the optimal value function defined by a polyhedral convex objective function and a polyhedral convex constrained multifunctions.

\medskip
The theorem below establishes a framework in which the composition of two generalized polyhedral convex multifunctions yields another generalized polyhedral convex multifunction.

\begin{Theorem}\label{composition}
Let $F\colon X \rightrightarrows Y$ and $G\colon Y \rightrightarrows Z$ be generalized polyhedral convex multifunctions whose graphs are given by
	\begin{eqnarray}\label{def_Graph_F}
		\begin{array}{ll}
			{\rm gph}\,F=\big\{(x, y) \in X \times Y \; \big|\; & A_1(x)+A_2(y)=u, \\
			&\langle x^*_i, x \rangle + \langle y^*_{1,i}, y \rangle \leq \alpha_i,\, i=1,\ldots,p\big\},
		\end{array}
	\end{eqnarray}
	\begin{eqnarray}\label{def_Graph_G}
		\begin{array}{ll}
			{\rm gph}\,G=\big\{(y, z) \in Y \times Z \; \big|\; &B_1(y)+B_2(z)=v, \\
			&\langle y^*_{2,j}, y \rangle + \langle z^*_j, z \rangle \leq \beta_j,\, j=1,\ldots,q\big\},
		\end{array}
	\end{eqnarray}		
where $A_1\colon X \to U$, $A_2\colon Y \to U$, $B_1\colon Y \to V$, $B_2\colon Z \to V$ are continuous linear mappings between locally convex Hausdorff topological vector spaces, $u \in U$, $v \in V$, $x^*_i \in X^*$, $y^*_{1,i} \in Y^*$, $\alpha_i \in \mathbb{R}$ for $i=1,\ldots,p$, $y^*_{2,j} \in Y^*$, $z^*_{j} \in Z^*$, $\beta_j \in \mathbb{R}$ for $j=1,\ldots,q$. If the continuous linear mapping $(A_2, B_1)\colon Y\to U\times V$ defined by $(A_2,B_1)(y)=(A_2(y),B_1(y))$ for all $y\in Y$ is closed under finite-codimensional  subspaces, then the multifunction $G \circ F\colon X \rightrightarrows Z$ is generalized polyhedral convex.
\end{Theorem}

\noindent
{\bf Proof.} Define the sets
\begin{align*}
	&\Omega_1=\big\{(x,z,y) \in X \times Z \times Y \; \big |\;  (x,y) \in {\rm gph}\,F \big\},\\
	&\Omega_2=\big\{(x,z,y) \in X \times Z \times Y \;\big | \;  (y,z) \in {\rm gph}\,G \big\}.
\end{align*}
We have
\begin{eqnarray*}
	\begin{array}{ll}
		\Omega_1 \cap \Omega_2=\big\{(x,z,y) \in X \times Z \times Y \; \big |\; &  A_1(x)+A_2(y)=u, B_1(y)+B_2(z)=v,\\
		&\langle x^*_i, x \rangle + \langle y^*_{1,i}, y \rangle \leq \alpha_i,\, i=1,\ldots,p,\\
		&\langle y^*_{2,j}, y \rangle + \langle z^*_j, z \rangle \leq \beta_j,\, j=1,\ldots,q\big\}.
	\end{array}
\end{eqnarray*}
Let $W= X \times Z$, $\widetilde{A}_1\colon W \to U \times V$, $\widetilde{A}_2\colon Y \to U \times V$, $\widetilde{u} \in W$, $\widetilde{x}^*_i, \widetilde{z}^*_j \in W^*$ for $i=1,\ldots,p$ and $j=1,\ldots,q$ be given by setting
\begin{equation*}
	\widetilde{A}_1 (w)= \widetilde{A}_1 \begin{pmatrix}
		x \\ z
	\end{pmatrix}=\begin{pmatrix}
		A_1(x) \\ B_2(z)
	\end{pmatrix}, \, \widetilde{A}_2(y)=\begin{pmatrix}
		A_2(y) \\ B_1(y)
	\end{pmatrix},
\end{equation*}
\begin{equation*}
	\widetilde{u}=\begin{pmatrix}
		u \\ v
	\end{pmatrix},\quad \langle \widetilde{x}^*_i, \begin{pmatrix}
		x \\ z
	\end{pmatrix} \rangle =  \langle x^*_i, x \rangle, \quad  \langle \widetilde{z}^*_j, \begin{pmatrix}
		x \\ z
	\end{pmatrix} \rangle =  \langle z^*_j, z \rangle
\end{equation*} for all $w=(x,z)\in X \times Z$ and $y\in Y$.
Then one has
\begin{eqnarray}\label{Intersection_Omega_1}
	\begin{array}{ll}
		\Omega_1 \cap \Omega_2=\big\{(w,y) \in W \times Y  \; \big|\; & \widetilde{A}_1(w)+\widetilde{A}_2(y)=\widetilde{u},\\
		&\langle \widetilde{x}^*_i, w \rangle + \langle y^*_{1,i}, y \rangle \leq \alpha_i,\, i=1,\ldots,p,\\
		&\langle y^*_{2,j}, y \rangle + \langle \widetilde{z}^*_j, w \rangle \leq \beta_j,\, j=1,\ldots,q\big\}.
	\end{array}
\end{eqnarray}
Let $\pi_{W}\colon \, W \times Y \to  W$ be the projection mapping from $W \times Y$ onto~$W$. Then it holds that $${\rm gph}(G \circ F)=\pi_{W} (\Omega_1 \cap \Omega_2).$$ In other words, ${\rm gph}(G \circ F)$ is the domain of the multifunction $T\colon W\rightrightarrows Y$ with
$$T(w)=\big\{y\mid (w,y)\in \Omega_1 \cap \Omega_2\big\},\ \; w\in W.$$ Since $\widetilde{A}_2=(A_2, B_1)$, the continuous linear mapping $\widetilde{A}_2$ is closed under finite-codimensional  subspaces by our assumptions. Therefore, by Theorem~\ref{Property_domF} and formula~\eqref{Intersection_Omega_1}  we can infer that $\pi_{W} (\Omega_1 \cap \Omega_2)$ is a generalized polyhedral convex set in~$W$. Therefore, $G \circ F$ is a generalized polyhedral convex multifunction.
$\hfill\Box$

\medskip
The theorem below states that the composition of two polyhedral convex multifunctions is itself a polyhedral convex multifunction.

\begin{Theorem}
If $F\colon X \rightrightarrows Y$ and $G\colon Y \rightrightarrows Z$ are polyhedral convex multifunctions, then $G \circ F\colon X \rightrightarrows Z$ is a polyhedral convex multifunction.
\end{Theorem}
\noindent
{\bf Proof.} We can assume that ${\rm gph}\,F$ and ${\rm gph}\,G$ are given by \eqref{def_Graph_F} and \eqref{def_Graph_G} with $U=V=\{0\}$, $A_1\equiv 0$, $A_2 \equiv 0$, $B_1\equiv 0$, $B_2 \equiv 0$ and $u=v=0$. Arguing similarly to the proof of Theorem~\ref{composition}, we obtain \eqref{Intersection_Omega_1} where $\widetilde{A}_1\equiv 0$ and $\widetilde{A}_2\equiv 0$. So, $\Omega_1 \cap \Omega_2$ is a polyhedral convex set in $W \times Y$. Then, applying Corollary~\ref{pcs_cr} for $\Omega_1 \cap \Omega_2$, one has $\pi_{W} (\Omega_1 \cap \Omega_2)$ is a polyhedral convex set in $W$. Since ${\rm gph}(G \circ F)=\pi_{W} (\Omega_1 \cap \Omega_2)$, we can assert that the multifunction $G \circ F\colon X \rightrightarrows Z$ is polyhedral convex.
$\hfill\Box$

\medskip
To conclude this section, we show that the sum of two polyhedral convex multifunctions is a polyhedral convex one and then follow up with an example showing that the conclusion no longer holds if the multifunctions involved are just generalized polyhedral convex.

\begin{Theorem}\label{sum_rule}
	If $F_1, F_2\colon  X \rightrightarrows Y$ are polyhedral convex multifunctions, then the multifunction $F_1+F_2$ is also polyhedral convex.
\end{Theorem}

\noindent
{\bf Proof.}
Consider the sets
\begin{eqnarray*}\label{O12}
	\begin{array}{ll}
		&\Omega_1=\big\{(x, y_1, y_2)\in X \times Y \times Y \; |\; y_1\in F_1(x)\big\} = ({\rm gph}\,F_1) \times Y,\\
		&\Omega_2=\{(x, y_1, y_2)\in X \times Y \times Y \; |\; y_2\in F_2(x)\}.
	\end{array}
\end{eqnarray*}
Since $F_1$ and $F_2$ are polyhedral convex multifunctions, $\Omega_1$ and $\Omega_2$ are polyhedral convex sets in $X \times Y \times Y$. Hence, $\Omega_1 \cap \Omega_2$ is a polyhedral convex set in $X \times Y \times Y$.

Let the continuous linear mapping $A\colon X \times Y \times Y \to X \times Y$ be given by
\begin{equation*}
	A(x,y_1,y_2) = (x,y_1+y_2)\ \; \mbox{\rm for}\ \; (x,y_1,y_2)\in X \times Y \times Y.
\end{equation*}
It is clear that ${\rm gph}(F_1+F_2) = A(\Omega_1\cap\Omega_2)$. If $\Omega_1 \cap \Omega_2=\emptyset$, then ${\rm gph}\,(F_1+F_2)=\emptyset$; so the multifunction $F_1+F_2$ is polyhedral convex. To proceed further, let us suppose that $\Omega_1 \cap \Omega_2$ is nonempty.  Since $\Omega_1 \cap \Omega_2$ is a polyhedral convex set in $X \times Y \times Y$, there exist $x^*_i \in X^*$, $y^*_{1,i} \in Y^*$, $y^*_{2,i} \in Y^*$, and $\alpha_i \in \mathbb{R}$ for $i=1,\ldots,m$ such that
\begin{equation*}
	\Omega_1 \cap \Omega_2=\big\{(x, y_1,y_2) \in X \times Y \times Y \;\big|\; \langle x^*_i, x \rangle + \langle y^*_{1,i}, y_1 \rangle + \langle y^*_{2,i}, y_2 \rangle \leq \alpha_i, \ i=1,\ldots,m\big\}.
\end{equation*}
Consider the sets
\begin{eqnarray*}
	\begin{array}{ll}
		&X_0=\{x \in X \, \mid \, \langle x^*_i, x \rangle =0,\, i=1,\ldots,m  \},\\
		&Y_{1,0}=\{y_1 \in Y \, \mid \, \langle y^*_{1,i}, y_1 \rangle=0,\, i=1,\ldots,m  \},\\
		&Y_{2,0}=\{y_2 \in Y \, \mid \, \langle y^*_{2,i}, y_2 \rangle=0, \, i=1,\ldots,m  \}.
	\end{array}
\end{eqnarray*}
Because $X_0 \subset X$, $Y_{1,0} \subset Y$, $Y_{2,0} \subset Y$ are closed linear subspaces of finite codimension, one can find finite-dimensional linear subspaces $X_1 \subset X$, $Y_{1,1} \subset Y$, and $Y_{2,1} \subset Y$ such that $$X=X_0 + X_1,\ \;  Y=Y_{1,0}+Y_{1,1},\ \; Y=Y_{2,0}+Y_{2,1},$$ $X_0 \cap X_1=\{0\}$, $Y_{1,0} \cap Y_{1,1}=\{0\}$, and $Y_{2,0} \cap Y_{2,1}=\{0\}$. According to~\cite[Theorem~1.21(b)]{Rudin_1991}, the subspaces $X_1, Y_{1,1}, Y_{2,1}$ are closed. It is clear that $$D_1=\big\{(x, y_1,y_2) \in X_1 \times Y_{1,1} \times Y_{2,1} \mid \langle x^*_i, x \rangle + \langle y^*_{1,i}, y_1 \rangle + \langle y^*_{2,i}, y_2 \rangle \leq \alpha_i, \ i=1,\ldots,m\big\}$$ is a polyhedral convex set in $X_1 \times Y_{1,1} \times Y_{2,1}$. Put $D_0=X_0 \times Y_{1,0} \times Y_{2,0}$. It is easy to verify that $$D_0+ D_1\subset \Omega_1 \cap \Omega_2.$$ The reverse inclusion is also true. Indeed, for each $(x,y_1, y_2) \in \Omega_1 \cap \Omega_2$ there exist $x_0 \in X_0$, $x_1 \in X_1$, $y_{1,0} \in Y_{1,0}$, $y_{1,1} \in Y_{1,1}$, $y_{2,0} \in Y_{2,0}$, $y_{2,1} \in Y_{2,1}$ satisfying $x=x_0+x_1$, $y_1=y_{1,0}+y_{1,1}$, $y_2=y_{2,0}+y_{2,1}$. Since
\begin{eqnarray*}
	\begin{array}{ll}
		\langle x^*_i, x_1 \rangle + \langle y^*_{1,i}, y_{1,1} \rangle + \langle y^*_{2,i}, y_{2,1} \rangle &= \big(\langle x^*_i, x \rangle - \langle x^*_i, x_0 \rangle \big) + \big(\langle y^*_{1,i}, y_1 \rangle - \langle y^*_{1,i}, y_{1,0} \rangle \big)\\
		&\qquad  + \big(\langle y^*_{2,i}, y_2 \rangle - \langle y^*_{2,i}, y_{2,0} \rangle \big)\\
		&=\langle x^*_i, x \rangle + \langle y^*_{1,i}, y_1 \rangle+ \langle y^*_{2,i}, y_2 \rangle \leq \alpha_i
	\end{array}	
\end{eqnarray*}
for every $i=1,\ldots,m$, it follows that $(x_1,y_{1,1},y_{2,1}) \in D_1$; so
$$(x,y_1,y_2)=(x_0,y_{1,0},y_{2,0}) + (x_1,y_{1,1},y_{2,1}) \in D_0+D_1.$$
We have thus proved that $\Omega_1 \cap \Omega_2=D_0+ D_1$. Hence
\begin{eqnarray*}
	\begin{array}{ll}
		A(\Omega_1\cap\Omega_2)&=A(D_0)+A(D_1)\\
		&=\big(X_0\times (Y_{1,0}+Y_{2,0})\big)+A(D_1).
	\end{array}	
\end{eqnarray*}
Since $D_1$ is a polyhedral convex set of the finite-dimensional space $X_1 \times Y_{1,1} \times Y_{2,1}$, invoking Theorem~19.1 in~\cite{Rockafellar_1970} one can represent $D_1$ as
\begin{eqnarray*}
	\begin{array}{ll}
		D_1=\Big\{ \sum\limits_{i=1}^k \lambda_i u_i + \sum\limits_{j=1}^\ell \mu_j v_j \; \big|\; & \lambda_i \geq 0, \ \forall i=1,\ldots,k,\sum\limits_{i=1}^k \lambda_i=1,\ \,    \\
		&\mu_j \geq 0,\ \forall j=1,\ldots,\ell \Big\},
	\end{array}
\end{eqnarray*}
where $u_i \in D_1$ for $i=1,\ldots,k$, and $v_j \in X_1 \times Y_{1,1} \times Y_{2,1}$ for $j=1,\ldots,\ell$. Then we have
\begin{eqnarray*}
	\begin{array}{ll}
		A(\Omega_1\cap\Omega_2)\!=\!\Big\{\!\sum\limits_{i=1}^k \lambda_i (A(u_i)) \!+\! \sum\limits_{j=1}^\ell \mu_j (A(v_j))\ \big| &  \lambda_i \geq 0, \, \forall i=1,\ldots,k,\ \sum\limits_{i=1}^k \lambda_i=1, \\
		& \mu_j \geq 0,\, \forall j=1,\ldots,\ell\Big\} \!  + \! \big(X_0\times (Y_{1,0}+Y_{2,0})\big).
	\end{array}
\end{eqnarray*}
Since $X_0 \subset X$, $Y_{1,0} \subset Y$, $Y_{2,0} \subset Y$ are closed linear subspaces of finite codimension, $X_0\times (Y_{1,0}+Y_{2,0})$ is a finite-codimensional closed linear subspace of $X \times Y \times Y$. Hence, Theorem~\ref{rep_gpcs} assures that $A(\Omega_1\cap\Omega_2)$ is a polyhedral convex set in~$X \times Y \times Y$. Since ${\rm gph}(F_1+F_2) = A(\Omega_1\cap\Omega_2)$, the set ${\rm gph}\,(F_1+F_2)$ is polyhedral convex. So, $F_1+F_2$ is a polyhedral convex multifunction.
$\hfill\Box$

\medskip
One may ask whether the statement in Theorem~\ref{sum_rule} applies to the summation of generalized polyhedral convex multifunctions. To clarify this, we now provide an example.

\begin{Example}{\rm Choose a suitable topological vector space~$X$ and closed linear subspaces $X_1, X_2$ of~$X$ satisfying $\overline{X_1+X_2}=X$ and $X_1+X_2 \neq X$ (see\cite[Remark~2.12]{Luan_Yao_Yen_2018} for more details). Let $F_1, F_2\colon X \rightrightarrows X$ be given by $F_1(x)=X_1, F_2(x)=X_2$ for all $x \in X$. It is clear that $F_1$ and $F_2$ are generalized polyhedral convex multifunctions. Since ${\rm gph}(F_1+F_2)=X \times (X_1+X_2)$ is not closed in the product space $X\times X$}, $F_1+F_2$ is not a generalized polyhedral convex multifunction.
\end{Example}

Given a function $\varphi\colon X \times Y \rightarrow \overline{\mathbb{R}}$ and a multifunction $F\colon X \rightrightarrows Y$, define the {\em optimal value function} $\mu\colon X\to \Bar{\mathbb{R}}$ associated with $\varphi$ and $F$ by
\begin{align}\label{marginalfunction}
	\mu(x)= \inf \big\{\varphi (x,y)\; \big|\; y \in F(x)\big\}, \ \; x\in X.
\end{align}
Here we use the convention $\inf \emptyset =\infty$. The \textit{solution map}  $M\colon X \rightrightarrows Y$ of the optimization problem in~\eqref{marginalfunction} is defined by
\begin{align}\label{solution_map}
	M(x)=\big\{y \in F(x)\; \big |\;  \mu(x)= \varphi (x,y)\big\}, \ \; x\in X.
\end{align}

The next result concerns the nonempty property of the solution set $M(x)$ at a given parameter $x\in X$.

\begin{Proposition}\label{Ex_Solution} Consider the optimal value function $\mu$ from \eqref{marginalfunction} and the solution mapping from \eqref{solution_map} in which $\varphi$ is a proper generalized polyhedral convex function and $F$ is a generalized polyhedral convex multifunction. For an element $x\in X$, if $\mu(x)$ is finite, then $M(x)$ is a nonempty subset of $Y$.
\end{Proposition}
{\bf Proof.} Fix $x\in X$ and assume that $\mu(x)$ is finite, i.e., $\mu(x)\in \mathbb R$. As $F$ is a generalized polyhedral convex multifunction, $F(x)$ is a generalized polyhedral convex set in $Y$ by Proposition~\ref{gpcs_image1}. Since $\gamma=\mu(x)$ is finite, we see that $F(x)$ is nonempty,  $\varphi(x,y) \ge \gamma$ for all $y \in F(x)$, and there exists $\overline{y} \in F(x)$ such that $\varphi(x,\overline{y})$ is finite.

Let the function $\varphi_x \colon Y \rightarrow \overline{\mathbb{R}}$ be given by $\varphi_x (y)=\varphi(x,y)$.  
Since $\varphi$ is a proper function and $\varphi(x,\overline{y})$ is finite, the function  $\varphi_x$ is also proper. Next, we will claim that $\varphi_x$ is a generalized polyhedral convex function. Since $\varphi$ is a proper generalized polyhedral convex function, the set $\mbox{\rm epi}\, \varphi$ can be represented by 
\begin{eqnarray*}\label{def_epi_varphi}
	\begin{array}{ll}
		\mbox{\rm epi}\, \varphi=\big\{(x, y,t) \in X \times Y \times \mathbb{R}\; \big |\; &B_1(x)+B_2(y)+B_3(t)=z,\\
		&\langle u^*_j, x \rangle + \langle v^*_j, y \rangle +\langle t_j, t \rangle \leq \alpha_j,\, j=1,\ldots,k\big\},
	\end{array}
\end{eqnarray*}	
where $B_1$ (resp., $B_2$, $B_3$) is a continuous linear mapping from $X$ (resp., from $Y$, from $\mathbb{R}$) to $Z$, $z \in Z$, $u^*_j \in X^*$, $v^*_j \in Y^*$, $t_j \in \mathbb{R}$, $\alpha_j \in \mathbb{R}$ for $j=1,\ldots,k$. Then one has
\begin{eqnarray*}
	\begin{array}{ll}
		\mbox{\rm epi}\, \varphi_x&=\big\{(y,t) \in Y \times \mathbb{R}\; \big |\; \varphi_x(y)\le t  \big\}\\
		&=\big\{(y,t) \in Y \times \mathbb{R}\; \big |\; \varphi(x,y)\le t  \big\}\\
		&=\big\{(y,t) \in Y \times \mathbb{R}\; \big |\; (x,y,t) \in \mbox{\rm epi}\, \varphi \big\}\\
		&=\big\{(y,t) \in Y \times \mathbb{R}\; \big |\; B_1(x)+B_2(y)+B_3(t)=z, \, \langle u^*_j, x \rangle + \langle v^*_j, y \rangle +\langle t_j, t \rangle \leq \alpha_j,\, j=1,\ldots,k\big\}\\
		&=\big\{(y,t) \in Y \times \mathbb{R}\; \big |\; B_2(y)+B_3(t)=z-B_1(x), \, \langle v^*_j, y \rangle +\langle t_j, t \rangle \leq \alpha_j-\langle u^*_j, x \rangle,\, j=1,\ldots,k\big\}
	\end{array}
\end{eqnarray*}	
It follows that $\varphi_x$ is a proper generalized polyhedral convex function and ${\rm dom}\,\varphi_x \cap F(x)$ is nonempty. Since $\varphi_x(y) \ge \gamma$ for all $y \in F(x)$, applying \cite[Theorem~3.1]{Luan_Yao_2019}, one can assert that the problem 
\begin{align*}
	&\mbox{\rm minimize }\;\varphi_x(y)\\
	&\mbox{\rm subject to }\;y\in F(x)
\end{align*}
has an optimal solution. Therefore, $M(x)$ is a nonempty set.
$\hfill\Box$

\medskip
The following proposition enables us to represent the epigraph of the optimal value function in terms of the image of a generalized polyhedral convex set under a projection mapping.

\begin{Proposition}\label{epi_mu} Consider the optimal value function $\mu$ from \eqref{marginalfunction} and let
	\begin{equation}\label{Omega1_2}
		\Omega_1=\mbox{\rm epi}\, \varphi\ \; \mbox{\rm and }\ \;\Omega_2=({\rm gph}\,F) \times \mathbb{R}.
	\end{equation}
	We have the representation
	\begin{equation}\label{relation}
		\overline{\mbox{\rm epi}\, \mu}=\overline{\pi_{X,\mathbb{R}}(\Omega_1 \cap \Omega_2)},
	\end{equation}
	where $\pi_{X,\mathbb{R}} \colon \, X \times Y \times \mathbb{R} \to X \times \mathbb{R}$ is the projection mapping from $X \times Y \times \mathbb{R}$ onto~$X \times \mathbb{R}$.
	If we assume in addition that $\varphi$ is a proper generalized polyhedral convex function and $F$ is a generalized polyhedral convex multifunction, then the closure signs in~\eqref{relation}
	can be omitted, i.e.,
	\begin{equation}\label{relation1}
		\mbox{\rm epi}\, \mu=\pi_{X,\mathbb{R}}(\Omega_1 \cap \Omega_2).
	\end{equation}
\end{Proposition}
\noindent
{\bf Proof.} Consider the set
\begin{equation*}
	\mbox{\rm epi}_s\, \mu=\big\{(x, \lambda)\in X\times \mathbb{R}\; \big |\; \mu(x)<\lambda\big\}.
\end{equation*}
We have
\begin{equation}\label{ic1}
	\mbox{\rm epi}_s \,\mu\subset \pi_{X,\mathbb{R}}(\Omega_1 \cap \Omega_2)\subset \mbox{\rm epi}\, \mu.
\end{equation}
Indeed, for any $(x, \lambda)\in \mbox{\rm epi}_s\,\mu$ we have $\mu(x)<\lambda$ and thus there exists $\bar{y}\in F(x)$ such that $\varphi(x, \bar{y})<\lambda$. Then we get $(x, \bar{y}, \lambda)\in \Omega_1\cap \Omega_2$, which implies that $(x, \lambda)\in \pi_{X,\mathbb{R}}(\Omega_1 \cap \Omega_2)$. This justifies the first inclusion in~\eqref{ic1}. To prove the second inclusion, take any $(x,\lambda) \in \pi_{X,\mathbb{R}}(\Omega_1 \cap \Omega_2)$. Then there is a point $\bar{y} \in Y$ such that $(x,\bar{y},\lambda)  \in \Omega_1 \cap \Omega_2$. It means that $\varphi(x,\bar{y}) \le \lambda$ and $\bar{y} \in F(x)$. Thus,  $\mu(x) \le \varphi(x,\bar{y})\le \lambda$. This implies that $(x,\lambda) \in \mbox{\rm epi}\, \mu$ and completes the proof of \eqref{ic1}. Finally, using \eqref{ic1} and the obvious equality $\overline{\mbox{\rm epi}_s\,\mu}=\overline{\mbox{\rm epi}\,\mu}$ gives us \eqref{relation}.

Now, assume that $\varphi$ is a proper generalized polyhedral convex function and $F$ is a generalized polyhedral convex multifunction. The proof above gives us
\begin{equation*}
	\pi_{X,\mathbb{R}}(\Omega_1 \cap \Omega_2)\subset\mbox{\rm epi}\, \mu.
\end{equation*}
Thus, to justifies \eqref{relation1}, it suffices to show that $\mbox{\rm epi}\, \mu\subset\pi_{X,\mathbb{R}}(\Omega_1 \cap \Omega_2)$. Take any $(x,\lambda) \in \mbox{\rm epi}\, \mu$ and get $\mu(x)\leq \lambda$. If $\mu(x)$ is finite, by Proposition~\ref{Ex_Solution}, there exists $y \in F(x)$ such that $\varphi(x,y)=\mu(x) \le \lambda$.  Now, consider the case where $\mu(x)=-\infty<\lambda$. In this case we can also choose $y\in F(x)$ such that $\varphi(x, y)<\lambda$. Thus, $(x,y,\lambda) \in \mbox{\rm epi}\, \varphi$ and $(x,y,\lambda) \in \Omega_2$. It follows that $(x,y,\lambda)  \in \Omega_1 \cap \Omega_2$. Therefore, $(x,\lambda) \in \pi_{X,\mathbb{R}}(\Omega_1 \cap \Omega_2)$, so~\eqref{relation1} is valid.  $\hfill\Box$

\medskip
The next theorem characterizes the generalized polyhedral convex property of the optimal value function~$\mu$ via its lower semicontinuity.

\begin{Theorem} Consider the optimal value function $\mu$ from \eqref{marginalfunction} in which $\varphi$ is a generalized polyhedral convex function and $F$ is a generalized polyhedral convex multifunction. The function $\mu$ is generalized polyhedral convex if and only if $\mu$ is lower semicontinuous on $X$.
\end{Theorem}
\noindent
{\bf Proof.} If $\mu$ is a generalized polyhedral convex function, then ${\rm epi}\,\mu$ is a generalized polyhedral convex set and hence ${\rm epi}\,\mu$ is closed. Thus, $\mu$ is lower semicontinuous on $X$.

Now, suppose that $\mu$ is lower semicontinuous on $X$. Then ${\rm epi}\,\mu$ is closed. Combining this fact with~\eqref{relation} gives us the equality
\begin{equation}\label{relation5}
	{\rm epi}\,\mu=\overline{\pi_{X,\mathbb{R}}(\Omega_1 \cap \Omega_2)},
\end{equation} where $\Omega_1$ and $\Omega_2$ are defined by~\eqref{Omega1_2}. Since $F$ is a generalized polyhedral convex multifunction and $\varphi$ is a generalized polyhedral convex function, the set $\Omega_1 \cap \Omega_2$ is generalized polyhedral convex. Applying Proposition~2.10 in \cite{Luan_Yao_Yen_2018}, we can conclude that  $\overline{\pi_{X,\mathbb{R}}(\Omega_1 \cap \Omega_2)}$ is a generalized polyhedral convex set. Therefore, by equality~\eqref{relation5} we can assert that $\mbox{\rm epi}\,\mu$ is generalized polyhedral, so $\mu$ is a generalized convex function.
$\hfill\Box$

\medskip
Sufficient conditions for the polyhedral convex property of the optimal value function $\mu$ are given in the following theorem.

\begin{Theorem}\label{mu_pcs} Consider the optimal value function $\mu$ from \eqref{marginalfunction}. If $\varphi$ is a proper polyhedral convex function and~$F$ is a polyhedral convex multifunction, then $\mu$ is a polyhedral convex function.
\end{Theorem}
\noindent
{\bf Proof.} The polyhedral convexity of $F$ and $\varphi$ guarantees that $\Omega_1 \cap \Omega_2$ is a polyhedral convex set in $X \times Y \times \mathbb{R}$. Using Corollary~\ref{pcs_cr},  we can assert that $\pi_{X,\mathbb{R}}(\Omega_1 \cap \Omega_2)$ is a polyhedral convex set in $X\times\mathbb{R}$. Combining this with~\eqref{relation1}, we conclude that $\mu$ is a polyhedral convex function.
$\hfill\Box$

\medskip
The conclusion of Theorem~\ref{mu_pcs} may not hold if one of the assumptions is violated. The next example shows that if $\varphi$ is a proper polyhedral convex function and $F$ is merely  a generalized polyhedral convex multifunction, then $\mu$ may not be a generalized polyhedral convex function.

\begin{Example}\label{Example2}{\rm Let $X, Y,$ and $F$ be as in Example~\ref{Example1}. Set $$X_1=\big\{x \in C[a,b] \; \big|\; x \text{ is continuously differentiable on } (a,b), \, x(a)=0 \big\}$$ and note that $X_1$ is a non-closed linear subspace of $X$. Since ${\rm gph}\,F=\Omega$, where $\Omega$ is defined by~\eqref{Omega_ex1}, we have $F(x)=\{-\dot{x}\}$ for all $x \in X_1$ with $\dot{x}$ denoting the Fr\'echet derivative of $x$, and $F(x)=\emptyset$ for all $x \notin X_1$. Consider the proper polyhedral convex function $\varphi$  with $\varphi(x,y)=0$ for all $(x,y) \in X \times Y$. As $\mu(x)=0$ for every $x \in X_1$ and $\mu(x)=\infty$ for any $x \notin X_1$, we see that $\mbox{\rm epi}\,\mu=X_1 \times [0,\infty)$. Since the latter set is non-closed, $\mu$ is not a generalized polyhedral convex function.	
}\end{Example}

\section{Generalized Relative Interiors of Generalized Polyhedral Convex Sets}\markboth{\centerline{\sc Generalized Relative Interiors of Generalized Polyhedral Convex Sets}}{\centerline{\sc N.~N.~Luan, N.~M.~Nam, and N.~D.~Yen}}
\setcounter{equation}{0}

 The notion of {\em relative interior} has been known to be useful for the study of convex analysis in finite dimensions. Its importance has motivated the development of new notions of generalized relative interiors in infinite dimensions. In this section, we show that several generalized relative interior concepts known in the literature do coincide for generalized polyhedral convex sets in infinite dimensions. We also obtain representations of such generalized relative interiors for the graphs of generalized polyhedral convex multifunctions.

Recall (see, e.g., ~\cite[Definition~2.168]{mordukhovich_nam_2021}) that the {\em relative interior}, the \textit{intrinsic relative interior}, and the \textit{quasi-relative interior} of a subset $\Omega$ of $X$ are defined respectively by
\begin{align*}
&\mbox{\rm ri}\, \Omega=\big\{a\in \Omega\; |\; \exists\; \mbox{\rm a neighborhood }\; V \; \mbox{\rm of the origin  such that}\; (a+V)\cap \overline{\mbox{\rm aff}\, \Omega}\subset \Omega\big\}.\\
&\mbox{\rm iri}\, \Omega=\big\{a\in \Omega\; |\; \mbox{\rm cone}(\Omega-a)\; \mbox{\rm is a linear subspace of }X\big\},\\
&\mbox{\rm qri}\, \Omega=\big\{a\in \Omega\; |\; \overline{\mbox{\rm cone}(\Omega-a)}\; \mbox{\rm is a linear subspace of }X\big\}.
\end{align*}

By~\cite[Theorem~2.174]{mordukhovich_nam_2021}, the following inclusions hold
\begin{equation}\label{inclusions}
	\mbox{\rm ri}\, \Omega\subset \mbox{\rm iri}\, \Omega\subset \mbox{\rm qri}\, \Omega.
\end{equation}

The theorem below shows that these generalized relative interior notions coincide for generalized polyhedral convex sets. It is a basis for obtaining the subsequent useful result about generalized polyhedral convex multifunctions.

\begin{Theorem}\label{Poly} Let $X$ be a locally convex Hausdorff topological vector space. Consider the generalized polyhedral convex set
	\begin{equation*}
		P=\big\{x\in X\; \big|\; \langle x^*_i, x\rangle \leq  \alpha_i\ \; \mbox{\rm for all }\; i=1, \ldots, m\big\}\cap L,
	\end{equation*}
	where $x^*_i\in X^*$, $\alpha_i\in \mathbb{R}$ for all $i=1, \ldots, m$, and $L$ is a closed affine subspace of $X$. Suppose that $P$ is nonempty. Then $\mbox{\rm ri}\, P$ is nonempty and we have the equalities
	\begin{equation}\label{formula_for_riP}
		\mbox{\rm qri}\, P=\mbox{\rm iri}\, P=\mbox{\rm ri}\, P=\big\{x\in P\; \big|\; \langle x^*_i, x\rangle <\alpha_i\  \mbox{\rm for all } i\in I\big\},
	\end{equation} where \begin{equation*}
	I=\big\{i=1, \ldots, m\; \big|\; \exists \hat{x}_i\in P\  \mbox{\rm such that }\langle x^*_i, \hat{x}_i\rangle <\alpha_i\big\}.
\end{equation*}
\end{Theorem}
{\bf Proof.} In the first part of the proof, we follow the proof of \cite[Proposition~2.197]{Bonnans_Shapiro_2000}, while providing more details.


First, consider the case where $I\neq\emptyset$. Fix an element $a\in P \subset L$. Denote by $N$ the unique linear subspace parallel to $\mbox{\rm aff}\, P$. Let us show that
\begin{equation}\label{prl}
	N=\big\{x\in X\; \big|\; \langle x^*_i, x\rangle=0 \ \, \mbox{\rm for all }\, i \in \{1, \ldots, m\}\setminus I\big\} \cap \big(L-a\big).
\end{equation}
Recall that  $N=\mbox{\rm cone}(P-P)=\mbox{\rm span}(P-a)$ and $a+N=\mbox{\rm aff}\, P$. One has
\begin{equation*}
	\langle x^*_i, a\rangle \leq  \alpha_i\ \; \mbox{\rm for all }\; i=1, \ldots, m.
\end{equation*}
Observe that $\langle x^*_i, a\rangle=\alpha_i$ if $i\notin I$. The set on the right-hand side of \eqref{prl}, which is denoted by $N_1$, is a closed linear subspace. For any $x\in P$ we have
\begin{equation*}
	\langle x^*_i, x\rangle=\alpha_i\ \, \mbox{\rm whenever }\, i\notin I,
\end{equation*}
which implies that
\begin{equation*}
	\langle x^*_i, x-a\rangle=\alpha_i-\alpha_i=0\ \, \mbox{\rm whenever }\, i\notin I.
\end{equation*}
In addition, it is clear that $x-a \in (P-a) \subset (L-a)$. Thus, $x-a\in N_1$ and hence $P-a\subset N_1$. Then $N=\mbox{\rm span}(P-a)\subset N_1$ because $N_1$ is a linear subspace.

To prove the reverse inclusion in \eqref{prl}, take any $x\in N_1$ and get $\langle x^*_i, x\rangle=0$ for all $i\notin I$. For every $i\in I$, choose $\hat{x}_i\in P$ such that $\langle x^*_i, \hat{x}_i\rangle <\alpha_i$. Denote by~$p$ be the number of elements of $I$ and define
\begin{equation*}
	\hat{x}=\dfrac{1}{p}\sum_{i\in I}\hat{x}_i.
\end{equation*} Then we have $\hat{x}\in P$ and $\langle x^*_i, \hat{x}\rangle <\alpha_i$ for all $i\in I$.
Fix any $j\in I$. Since $\hat{x}_i\in P$, we see that $\langle x^*_j, \hat{x}_i\rangle \leq  \alpha_j$ if $i\neq j$, and $\langle x^*_j, \hat{x}_j\rangle <\alpha_j$. Therefore,
\begin{equation*}
	\langle x^*_j, \hat{x}\rangle=\dfrac{1}{p} \sum_{i\in I} \langle x^*_j, \hat{x}_i\rangle < \dfrac{1}{p} p \alpha_j=\alpha_j.
\end{equation*}
We have thus shown that $\langle x^*_j, \hat{x}\rangle <\alpha_j$ for all $j\in I$. Then, for a sufficiently small $t>0$, we have
\begin{equation*}
	\langle x^*_i, \hat{x}+tx\rangle <\alpha_i\ \; \mbox{\rm for all }i\in I,\ \mbox{\rm and }\; \langle x^*_i, \hat{x}+tx\rangle =\langle x^*_i, \hat{x}\rangle \leq  \alpha_i\ \; \mbox{\rm for all }\; i\notin I.
\end{equation*}
In addition, since $x \in L-a=L-\hat{x}$, one has $\hat{x}+x \in L$. Hence, $\hat{x}+tx=(1-t)\hat{x}+t(\hat{x}+x) \in L$ as $L$ is an affine subspace. Thus, $\hat{x}+tx\in P$; so $x\in \frac{1}{t}(P-\hat{x})\subset \mbox{\rm cone}(P-P)=N$, which completes the proof of~\eqref{prl}.

For convenience, let
\begin{equation}\label{C}
	C=\big\{x\in P\; |\; \langle x^*_i, x\rangle <\alpha_i\ \,\mbox{\rm for all }\, i\in I\big\}.
\end{equation}
We will show that $C=\mbox{\rm ri}\, P$. Taking any $x_0\in C$, we have $x_0\in P$ and
\begin{equation*}
	\langle x^*_i, x_0\rangle <\alpha_i\ \; \mbox{\rm for all }\; i\in I.
\end{equation*}
By the continuity of $x^*_i$ for $i\in I$, we can find a neighborhood $U$ of the origin such that
\begin{equation}\label{neighborhood}
	\langle x^*_i, x_0+u\rangle \leq  \alpha_i\ \, \mbox{\rm for all }\, u\in U \ \; \mbox{\rm and for all }\, i\in I.
\end{equation}
Let us show that
\begin{equation}\label{inclusion_P}
	(x_0+U)\cap \overline{\mbox{\rm aff}\, P}=(x_0+U)\cap (\hat{x}+N)\subset P.
\end{equation}
Note that $\hat{x}$ is chosen above and $N$ is closed with $\hat{x}+N=\mbox{\rm aff}\, P=\overline{\mbox{\rm aff}\, P} \subset L$ by the definition of parallel subspace. Hence, the equality in~\eqref{inclusion_P} is valid. To obtain the inclusion in~\eqref{inclusion_P}, fix any $x\in (x_0+U)\cap (\hat{x}+N)$. Then $x=x_0+u$ for some $u\in U$, and $x=\hat{x}+v \in L$ for some $v\in N$. For $i\in I$, from~\eqref{neighborhood} it follows that
\begin{equation*}
	\langle x^*_i, x\rangle=\langle x^*_i, x_0+u\rangle \leq  \alpha_i.
\end{equation*}
If $i\notin I$, by~\eqref{prl} we have
\begin{equation*}
	\langle x^*_i, x\rangle=\langle x^*_i, \hat{x}+v\rangle=\langle x^*_i, \hat{x}\rangle +\langle x^*_i, v\rangle=\langle x^*_i, \hat{x}\rangle \leq  \alpha_i.
\end{equation*}
It follows that $x\in P$, and so~\eqref{inclusion_P} holds. Then we get $x_0\in \mbox{\rm ri}\, P$, which justifies the inclusion $C\subset \mbox{\rm ri}\, P$.

Now we show that $\mbox{\rm ri}\, P\subset C$. Fix any $x_0\in \mbox{\rm ri}\, P$ and find a neighborhood $U$ of the origin such that
\begin{equation}\label{ri}
	(x_0+U)\cap \overline{\mbox{\rm aff}\, P}\subset P.
\end{equation}
By contradiction, suppose that $x_0\notin C$, and so there exists $j\in I$ such that $\langle x^*_j, x_0\rangle\geq \alpha_j$, which implies $\langle x^*_j, x_0\rangle=\alpha_j$. Since $U$ is a neighborhood of the origin, we can find $t>0$ sufficiently small such that $$z=x_0+t(x_0-\hat{x}_j)\in x_0+U,$$ where $\langle x^*_j, \hat{x}_j\rangle <\alpha_j$. Obviously, $z\in \overline{\mbox{\rm aff}\, P}$ because $z=-t\hat{x}_j+(1+t)x_0$, $\hat x_j\in P$, and $x_0\in P$. So, by~\eqref{ri} one has $z \in P$. This implies that $\langle x^*_j, z\rangle \le \alpha_j$. Then $x_0=\frac{1}{1+t}z+\frac{t}{1+t}\hat{x}_j$ and thus
\begin{equation*}
	\alpha_j= \langle x^*_j , x_0\rangle =\frac{1}{1+t} \langle x^*_j, z\rangle +\frac{t}{1+t}\langle x^*_j, \hat{x}_j\rangle <\frac{1}{1+t} \alpha_j+\frac{t}{1+t}\alpha_j=\alpha_j,
\end{equation*}
which is a contradiction. Therefore, $x_0\in C$, and so $\mbox{\rm ri}\, P\subset C$. We have thus proved that if $I\neq\emptyset$, then  $\mbox{\rm ri}\, P= C$.

Now, consider the case where $I=\emptyset$. In this case, we have
\begin{equation*}
	P=\big\{x\in X\; \big|\; \langle x^*_i, x\rangle =  \alpha_i\ \; \mbox{\rm for all }\; i=1, \ldots, m\big\}\cap L.
\end{equation*} It follows that $P=\mbox{\rm aff}\, P=\overline{\mbox{\rm aff}\, P}$. Therefore,  $\mbox{\rm ri}\, P=P$.
On the other hand, by~\eqref{C} we get $C=P$. Thus, the equality $\mbox{\rm ri}\, P= C$ is also valid in the case where $I\neq\emptyset$.

The preceding proof shows that $\mbox{\rm ri}\, P\neq\emptyset$.

By~\eqref{inclusions}, to obtain~\eqref{formula_for_riP}, it suffices to show that $\mbox{\rm qri}\, P\subset C=\mbox{\rm ri}\, P$. If $I=\emptyset$, then $C=\mbox{\rm ri}\, P=P$; hence the latter is valid. Now, consider the case where $I\neq\emptyset$ and suppose on the contrary that there is $a\in \mbox{\rm qri}\, P$ but $a\notin C$. Then, by~\eqref{C}, there exists $j\in I$ such that $\langle x^*_j, a\rangle=\alpha_j$. Choose $\hat{x}_j\in P$ such that $\langle x^*_j, \hat{x}_j\rangle <\alpha_j$. Obviously,
\begin{equation*}
	\hat{x}_j-a\in \overline{\mbox{\rm cone}(P-a)}.
\end{equation*}
Since $\overline{\mbox{\rm cone}(P-a)}$ is a linear subspace, we see that
\begin{equation*}
	a-\hat{x}_j\in \overline{\mbox{\rm cone}(P-a)}
\end{equation*}
For any $x\in P$, we have $\langle x^*_j, x-a\rangle=\langle x^*_j, x\rangle -\langle x^*_j, a\rangle \leq  \alpha_j-\alpha_j=0$, and hence $\langle x_j^*, z\rangle \leq  0$ for all $z\in \mbox{\rm cone}(P-a)$. By the continuity of $x^*_j$, we deduce that $\langle x^*_j, z\rangle \leq 0$ for all $z\in \overline{\mbox{\rm cone}(P-a)}$. Then
\begin{equation*}
	\langle x^*_j, a-\hat{x}_j\rangle \leq  0,
\end{equation*}
which  yields $\alpha_j=\langle x^*_j, a\rangle \leq \langle x^*_j, \hat{x}_j\rangle <\alpha_j$, a contradiction. This completes the proof. $\hfill\square$
\begin{Remark} {\rm The fact that the inequalities $\mbox{\rm ri}\, P=\mbox{\rm iri}\, P=\mbox{\rm qri}\, P$  hold for any generalized polyhedral convex set follows from the second assertion of Theorem~2.174 in~\cite{mordukhovich_nam_2021} and Proposition~2.197 from~\cite{Bonnans_Shapiro_2000}.}
\end{Remark}

To continue, we recall the following important properties of $\mbox{\rm iri}\, \Omega$ and $\mbox{\rm qri}\,\Omega$ for a convex set~$\Omega$ (see \cite[Propositions~2.169 and~2.181]{mordukhovich_nam_2021} more details). Given $\bar{x}\in \Omega$, we have
\begin{equation*}
\begin{aligned}
	& [\bar{x}\in \mbox{\rm iri}\, \Omega] \Longleftrightarrow [\forall x\in \Omega, \exists x^\prime\in \Omega\; \mbox{\rm such that }\bar{x}\in (x, x^\prime)],\\
	& [\bar{x}\notin \mbox{\rm qri}\,\Omega] \Longleftrightarrow [\{\bar{x}\}\; \mbox{\rm and }\Omega\; \mbox{\rm can be properly separated}].
\end{aligned}
\end{equation*}

The next theorem allows us to obtain a representation of the relative interior of a generalized polyhedral convex multifunction. This representation is based on Theorem~\ref{Poly} and the idea for proving Theorem~4.3 in~\cite{chuong-qri}. The closed range assumption is essential for the validity of the conclusion.

\begin{Theorem} \label{ri gph_thm} If the graph of a generalized polyhedral convex multifunction $F\colon X \rightrightarrows Y$ is described by~\eqref{def_Graph_F_1} in which the mapping $A_2$ is closed under finite-codimensional subspaces, then
	\begin{equation}\label{ri gph}
		\mbox{\rm ri}(\mbox{\rm gph}\, F)=\big\{(x, y)\; \big|\; x\in \mbox{\rm ri}(\mbox{\rm dom}\, F), \; y\in \mbox{\rm ri}(F(x))\big\}.
	\end{equation}
\end{Theorem}
{\bf Proof.} Suppose that the graph of $F$ is described by~\eqref{def_Graph_F_1} in which the mapping $A_2$ is closed under finite-codimensional subspaces. Then $\mbox{\rm dom}\, F$ is a generalized polyhedral convex set by Theorem~\ref{Property_domF}. Now, take any $(x_0, y_0)\in \mbox{\rm ri}(\mbox{\rm gph}\, F)$. First, let us show that $x_0\in \mbox{\rm ri}(\mbox{\rm dom}\, F)$. For any $x\in \mbox{\rm dom}\, F$ we can choose $y\in F(x)$, so $(x, y)\in \mbox{\rm gph}\, F$. Since $\mbox{\rm ri}(\mbox{\rm gph}\, F)=\mbox{\rm iri} (\mbox{\rm gph}\, F)$ by Theorem~\ref{Poly}, we can choose $(x^\prime, y^\prime)\in \mbox{\rm gph}\, F$ and $t\in (0,1)$ such that
\begin{equation*}
	(x_0, y_0)=t(x, y)+(1-t)(x^\prime, y^\prime).
\end{equation*}
Then $x_0\in (x, x^\prime)$, where $x^\prime\in \mbox{\rm dom}\, F$. Since $x\in \mbox{\rm dom}\, F$ can be chosen arbitrarily, it follows that $$x_0\in  \mbox{\rm iri}\, (\mbox{\rm dom}\, F).$$ So, applying Theorem~\ref{Property_domF} to the generalized polyhedral convex set $\mbox{\rm dom}\, F$ yields $x_0\in \mbox{\rm ri}(\mbox{\rm dom}\, F)$.  Now, let us show that $y_0\in \mbox{\rm ri}(F(x_0))$.  Observe that $F(x_0)=\mbox{\rm gph}\, F \cap \big(\{x_0\} \times Y\big)$ is a generalized polyhedral convex set. Take any $y\in F(x_0)$ and get $(x_0, y)\in \mbox{\rm gph}\, F$ and thus we can find $(x_1, y_1)\in \mbox{\rm gph}\, F$ and $s\in (0,1)$ such that
\begin{equation*}
	(x_0, y_0)=s (x_0, y)+(1-s)(x_1, y_1).
\end{equation*}
Then $x_1=x_0$ and $y_0\in (y, y_1)$, where $y_1\in F(x_0)$. Thus, $y_0\in \mbox{\rm iri}(F(x_0))=\mbox{\rm ri}(F(x_0))$. This justifies the inclusion $\subset$ in~\eqref{ri gph}.

We will now prove the inclusion $\supset$ in \eqref{ri gph}. Take any $x_0\in \mbox{\rm ri}(\mbox{\rm dom}\, F)$ and $y_0\in \mbox{\rm ri}(F(x_0))$. By contradiction, suppose that $(x_0, y_0)\notin \mbox{\rm ri}(\mbox{\rm gph}\, F)=\mbox{\rm qri}(\mbox{\rm gph}\, F)$, where the last equality holds by Theorem~\ref{Property_domF}.  By the proper separation mentioned prior to the formulation of this theorem, there exist $x^*\in X^*$ and $y^*\in Y^*$ such that
\begin{equation}\label{EQ1}
	\langle x^*, x\rangle +\langle y^*, y\rangle \leq  \langle x^*, x_0\rangle +\langle y^*, y_0\rangle
\end{equation}
for all $(x,y) \in \mbox{\rm gph}\, F$ and there exists $(\hat{x}, \hat{y})\in \mbox{\rm gph}\, F$ such that
\begin{equation}\label{EQ2}
	\langle x^*, \hat x\rangle +\langle y^*, \hat y\rangle < \langle x^*, x_0\rangle +\langle y^*, y_0\rangle .
\end{equation}
Substituting $(x,y)=(x_0,y)$, where $y\in F(x_0)$, to~\eqref{EQ1} gives us $\langle y^*, y\rangle \leq  \langle y^*, y_0\rangle$ for all $y\in F(x_0)$. Since $\hat x\in \mbox{\rm dom}\, F$ and $x_0\in \mbox{\rm ri}(\mbox{\rm dom}\, F)=\mbox{\rm iri}(\mbox{\rm dom}\, F)$, we can find $x_2\in \mbox{\rm dom}\, F$  and $\lambda\in (0,1)$ such that $x_0=\lambda\hat{x}+(1-\lambda)x_2$. Choosing $y_2\in F(x_2)$ and letting $y^\prime=\lambda \hat y+(1-\lambda)y_2$ give us
\begin{equation*}
	(x_0, y^\prime)=\lambda (\hat x, \hat y)+(1-\lambda )(x_2, y_2)\in \mbox{\rm gph}\, F
\end{equation*}
due to the convexity of $\mbox{\rm gph}\, F$, so $y^\prime\in F(x_0)$. Using \eqref{EQ1}, we have
\begin{equation}\label{EQ3}
	\langle x^*, x_2\rangle +\langle x^*, y_2\rangle \leq  \langle x^*, x_0\rangle +\langle y^*, y_0\rangle .
\end{equation}
Multiplying both sides of \eqref{EQ2} with $\lambda$, multiplying \eqref{EQ3} with $(1-\lambda)$, and adding the resulting inequalities give us
\begin{equation*}
	\langle x^*, \lambda\hat{x}+(1-\lambda)x_2\rangle +\langle y^*, \lambda \hat y+(1-\lambda)y_2\rangle <\langle x^*, x_0\rangle +\langle y^*, y_0\rangle .
\end{equation*}
Then we get
\begin{equation*}
	\langle x^*, x_0\rangle +\langle y^*, y^\prime\rangle <\langle x^*, x_0\rangle +\langle y^*, y_0\rangle,
\end{equation*}
which implies that $\langle y^*, y^\prime\rangle <\langle y^*, y_0\rangle$, where $y^\prime\in F(x_0)$. Remembering that $\langle y^*, y\rangle \leq  \langle y^*, y_0\rangle$ for all $y\in F(x_0)$, we see that $\{y_0\}$ and $F(x_0)$ can be properly separated, so $y_0\notin \mbox{\rm qri}(F(x_0))=\mbox{\rm ri}(F(x_0))$, which is a contradiction. This completes the proof. $\hfill\square$

\medskip
Thanks to Theorem~\ref{Property_domF} and Theorem~\ref{Poly}, we can obtain the following representations for the quasi-relative interior and intrinsic relative interior of the graph of a generalized polyhedral convex multifunction as a direct consequence of Theorems~\ref{ri gph}.

\begin{Corollary} \label{ri gph1} If the graph of $F$ is described by~\eqref{def_Graph_F_1} in which the mapping $A_2$ is closed under finite-codimensional subspaces, then
	\begin{align*}
		&\mbox{\rm qri}(\mbox{\rm gph}\, F)=\{(x, y)\; |\; x\in \mbox{\rm qri}(\mbox{\rm dom}\, F), \; y\in \mbox{\rm qri}(F(x))\},\\
		&\mbox{\rm iri}(\mbox{\rm gph}\, F)=\{(x, y)\; |\; x\in \mbox{\rm iri}(\mbox{\rm dom}\, F), \; y\in \mbox{\rm iri}(F(x))\}.
	\end{align*}
\end{Corollary}

Regarding Theorem~\ref{ri gph_thm}, an interesting question arises: {\textit{Can the assumption that the mapping $A_2$ is closed under finite-codimensional subspaces be removed from the theorem?} In order to answer this question in the negative, let us consider an example.

\begin{Example}\label{Example3}{\rm Let the spaces $X, Y$ and the multifunction $F$ be as in Example~\ref{Example1}. Since $\mbox{\rm gph}\, F$ is a closed linear subspace of $X \times Y$, one has $\mbox{\rm ri}(\mbox{\rm gph}\, F)=\mbox{\rm gph}\, F$. Here $$\mbox{\rm dom}\, F=\big\{ x \in C[a,b] \; \big | \; x \text{ is continuously differentiable on } (a,b),\, x(a)=0 \big\}$$ is a non-closed linear subspace, which is dense in $X$ (see \cite[Example 2.1]{Luan_APAN_2019} for details). Hence $$\mbox{\rm ri}(\mbox{\rm dom}\, F)=\mbox{\rm int}(\mbox{\rm dom}\, F)=\emptyset.$$  Consequently, the equality~\eqref{ri gph} does hold for the generalized polyhedral convex multifunction $F$ under our consideration.}
\end{Example}

\begin{acknowledgements}
Nguyen Ngoc Luan was funded by the Postdoctoral Scholarship Programme of Vingroup Innovation Foundation (VINIF) code VINIF.2022.STS.41. Nguyen Mau Nam would like to thank the Vietnam Institute for Advanced Study in Mathematics (VIASM) for hospitality. Valuable suggestions of the two anonymous referees are gratefully acknowledged. Among other things, our consideration of continuous linear mappings closed under finite-codimensional subspaces has the origin in an insightful discussion of one referee on the original proof of Theorem~\ref{Property_domF}.
\end{acknowledgements}

\end{document}